\documentclass[12pt]{article}

\usepackage{array}
\usepackage{graphicx}
\usepackage{amsmath}
\usepackage{hhline}

\newcommand{\be}{\begin{equation}}
\newcommand{\ee}{\end{equation}}
\newcommand{\bsube}{\begin{subequations}}
\newcommand{\esube}{\end{subequations}}
\newcommand{\ba}{\begin{array}}
\newcommand{\ea}{\end{array}}
\newcommand{\To}{\rightarrow}
\newcommand{\bea}{\begin{eqnarray}}
\newcommand{\eea}{\end{eqnarray}}
\newcommand{\bc}{\begin{center}}
\newcommand{\ec}{\end{center}}

\newcommand{\const}{{\rm const}}
\newcommand{\tu}{\tilde{u}}
\newcommand{\tv}{\tilde{v}}
\newcommand{\tw}{\tilde{w}}

\newcommand{\dz}{\Delta z}
\newcommand{\F}{{\mathcal F}}
\newcommand{\opi}{\omega_{\pi}}
\newcommand{\dopi}{\delta \omega_{\pi}}
\newcommand{\sech}{{\rm sech}\,}
\newcommand{\usol}{u_{\rm sol}}
\newcommand{\Lamr}{\Lambda_{R}}
\newcommand{\Lami}{\Lambda_{I}}

\renewcommand{\theequation}{\thesection.\arabic{equation}}

\title{\bf Stability analysis of the split-step Fourier method on the background of
a soliton of the nonlinear Schr\"odinger equation}

\author{ T.I. Lakoba\thanks{Department of Mathematics and Statistics, 16 Colchester Ave.,
 University of Vermont, Burlington, VT 05401, USA ({\tt lakobati@cems.uvm.edu}).  }  }

\begin{document}

\maketitle

\begin{abstract}
We analyze a numerical instability that occurs in the well-known split-step Fourier method
on the background of a soliton. This instability is found to be
very sensitive to small changes
of the parameters of both the numerical grid and the soliton, unlike the instability 
of most finite-difference schemes.
%  on the background of a monochromatic wave, considered earlier in the literature. 
Moreover,
the principle of ``frozen coefficients", in which variable coefficients are treated as
``locally constant" for the purpose of stability analysis, is strongly violated for
the instability of the split-step method on the soliton background. Our analysis 
explains all these features. It is enabled by the fact that the period of 
oscillations of the unstable Fourier modes is much smaller than the width of the soliton.
\end{abstract}

\bigskip


\bigskip


\newpage

% -------------------------------------------------------------------------
\section{Introduction}

The split-step Fourier method is widely used in numerical simulations of nonlinear wave
equations. In particular, it is the mainstream method in nonlinear optics, where the fundamental
equation describing propagation of an electromagnetic pulse or beam is the nonlinear Schr\"odinger equation
\be
iu_z  - \beta u_{tt} + \gamma u|u|^2 = 0\,, \qquad u(t,0)=u_0(t).
\label{e1_01}
\ee
In this paper, we use notations adopted in fiber optics \cite{Agrawal_book}, whereby $u$ is proportional
to the complex envelope of the electric field, $z$ is the propagation distance along the fiber, and 
$t$ is the time in the reference frame moving with the pulse. (The subscripts denote partial differentiation.)
Thus, $z$ and $t$ are the evolution and spatial variables, respectively.
Also, in Eq.~\eqref{e1_01}, $\beta$ and $\gamma$ are proportional to the group velocity dispersion 
and nonlinear refractive index of the optical fiber, respectively. Although these real-valued
constants can be scaled out of the
equation by an appropriate nondimensionalization, we will keep them in our analysis to distinguish 
contributions of the dispersive and nonlinear terms. Solitons with \ $u(z,|t|\To\infty)\To 0$ \ 
exist in \eqref{e1_01} for \ $\beta\gamma<0$.

The idea of the split-step method is the following.
Equation \eqref{e1_01} can be solved analytically when only one of the dispersive and nonlinear terms is
nonzero. Then, the approximate numerical solution of \eqref{e1_01} can be obtained in a sequence of
steps which alternatingly account either for the dispersion or for the nonlinearity:
\be
\ba{ll}
  \mbox{for $n$ from $1$ to $n_{\max}$ do:} & \\
    & \hspace*{-4cm} u_{\rm nonlin}(t) = u_n(t)\,\exp\big(i\gamma |u_n(t)|^2 \dz \big) \\
    & \hspace*{-4cm} u_{n+1} = \F^{-1}\Big[ \exp(i\beta\omega^2\; \dz) \, 
      \F\big[u_{\rm nonlin}(t) \big] \Big] \\
  \mbox{end} & 
\ea
\label{e1_02}
\ee
Here $\dz$ is the step size along the evolution variable, and the subscript $n$ now denotes
the value at the $n$th step\footnote{Using subscripts to denote iteration steps, as in \eqref{e1_02},
and partial differentiation, as in \eqref{e1_01}, will not lead to confusion.}.
In algorithm \eqref{e1_02}, the evolution due to the dispersive term is computed using the Fourier
transform $\F$ and its inverse $\F^{-1}$, where, e.g.:
\be
\F[u](\omega) =  \int_{-\infty}^{\infty} u(t) e^{-i\omega t} dt\,.
\label{e1_03}
\ee
For equations of the form \eqref{e1_01}, such an algorithm was originally applied in \cite{HardinTappert_73} --
\cite{FisherBischel_73} and later comprehensively studied in \cite{AblowitzTaha_84}.
Implementations where the dispersive term is computed using finite-difference discretization in $t$
can be used as well, but the Fourier-based implementation is preferred in optics, not only because of its
higher accuracy (exponential versus algebraic in $\Delta t$ for smooth pulses), but also because
it allows one to easily handle more complicated dispersive terms than in \eqref{e1_01}. In this paper,
we specifically focus on the Fourier-based implementation \eqref{e1_02} of the split-step method.

Since the split-step method is explicit, it can only be conditionally stable. 
However, the standard von Neumann stability analysis on the background of the trivial solution
$u(t)=0$ does not reveal any instability.
Weideman and Herbst \cite{WH} were the first to rigorously show that the split-step method is indeed
conditionally stable, with the instability being able to develop over the background of a {\em finite-amplitude}
monochromatic wave
\be
u_{\rm cw}=A\,\exp(iKz-i\Omega_{\rm cw} t), \qquad 
K=\beta\Omega_{\rm cw}^2 + \gamma |A|^2, \quad A=\const.
\label{e1_04}
\ee
A few years later, Matera et al \cite{Matera_93} obtained a similar result by a more heuristic method.

\begin{figure}[h]
 \vspace{-1cm}
\rotatebox{0}{\resizebox{7cm}{9cm}{\includegraphics[0in,0.5in]
 [8in,10.5in]{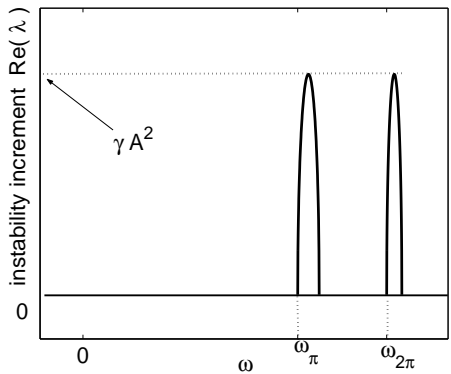}}}
\vspace{-2cm}
\caption{Schematics of the numerical
instability increment of method \eqref{e1_02} over the background of a monocromatic
wave \eqref{e1_04} as a function of frequency $\omega$. Only the part for $\omega>0$ is
shown, since the graph is symmetric about $\omega=0$. The width of the instability peak near
$\omega_{k\pi}$, where $k=1,2,\ldots$, equals \ $|\gamma A^2/(\beta \omega_{k\pi})|$ 
(see Appendix 3). Note that the case corresponding to
$\beta<0$ is shown; in the case of $\beta>0$, the peaks will occur on the opposite
sides of $\opi$ and $\omega_{2\pi}$.
}
\label{fig_1}
\end{figure}

The numerical
instability increment of the split-step method of Eq.~\eqref{e1_01} with the background solution 
$\bar{u}\equiv u_{\rm cw}$ \eqref{e1_04}, 
found by Weideman and Herbst, is schematically shown in Fig.~\ref{fig_1} for a
particular case of $\Omega_{\rm cw}=0$. (The formulae describing the location and growth rate
of unstable modes are given in Appendix 3.)
%   \footnote{
%  While the result shown in this figure can be deduced from formula (65) of 
%  \cite{WH}, we will later obtain it with a different method.}. 
Here the notations are the following. The numerical solution
of \eqref{e1_01} is assumed to have the form
\bsube
\be
u_n=\bar{u}+\tu_n, \qquad |\tu_n| \ll |\bar{u}|,
\label{e1_05a}
\ee
\be
\tu_n=\tilde{A}\exp(\lambda z_n - i\omega t), \qquad z_n=n\dz, \quad \omega=\frac{2\pi \,\ell}{T},
\label{e1_05b}
\ee
\label{e1_05}
\esube
where \ $-T/2 \le t \le T/2$ and the limits for the integer index $\ell$ are determined by the
number of grid points. In Fig.~\ref{fig_1}, 
\be
\opi=\sqrt{\frac{\pi}{|\beta|\,\dz}}, \qquad \mbox{i.e.} \quad |\beta|\opi^2\dz=\pi,
\label{e1_06}
\ee
and, similarly, $|\beta|\omega_{2\pi}^2\dz=2\pi$. Instability peaks of the same height
exist also near $\pm \omega_{3\pi}$, as long as  $\omega_{\max}\equiv \pi/\Delta t$,
where $\Delta t$ is the mesh size along $t$, extends beyond $\omega_{3\pi}$; etc.

Note that the instability depicted in Fig.~\ref{fig_1} is mild, in the sense that 
the numerical error grows by a factor $\exp\big(O(1)\big)$ over propagation distances 
of order $O(1)$. 
%  This is similar to the instability of the leap-frog numerical scheme
%  applied to the Heat equation when the step size along the evolution variable is
%  sufficiently small (see, e.g., \cite{NumMethods_somebook}). 
Furthermore, this instability remains mild even when $\dz$ significantly exceeds
its threshold
\be
\dz_{\rm thresh}\approx \frac{\pi}{|\beta| \omega_{\max}^2}
\label{e1_07}
\ee
(see \eqref{e1_06}). 
This should be contrasted with the behavior of finite-difference explicit schemes,
e.g., the Runge--Kutta scheme for the Heat equation, 
where the instability becomes strong (i.e. the numerical error grows by a factor of $\exp\big(O(1)\big)$
over propagation distances of only $O(\dz)$)
whenever the step size along the evolution variable exceeds the instability threshold.
%
%  In contrast, when the step size along the evolution variable
%  exceeds the corresponding threshold, the instability in the leap-frog scheme
%  (and other explicit schemes, such as the explicit Euler method etc.)
%  becomes strong, i.e. the numerical error grows by a factor of $\exp\big(O(1)\big)$
%  when the evolution variable increases by only $O(\dz)$. 
Another remarkable property of
the instability in Fig.~\ref{fig_1}
is that it occurs only in a finite (and rather narrow, of width $O(\gamma A^2/(|\beta|\opi)\,=\,O(\sqrt{\dz})$) 
band of spatial frequencies $\omega$,
whereas instabilities of other explicit schemes occur, typically, for {\em all} frequencies higher
than a certain minimum value. 
Let us note a curious consequence of this spectral selectivity of the instability of the
split-step method:
The method can be found stable even when $\dz$ exceeds the
threshold \eqref{e1_07}. This can occur when the width of the instability peak is
less than the mesh size $\Delta\omega$ in the frequency domain \cite{WH}; see also Appendix 3.

In this paper, we will show, first via numerical simulations and then by an analytical
calculation, that the properties of the instability of 
the split-step Fourier method on the background of a soliton solution are considerably 
different from properties of such instability on the monochromatic wave background.
In particular, these new properties are even more distinct from some propeties of 
instabilities of well-known finite-difference explicit schemes. Highlights of these new
properties include: \ high sensitivity to (i) the step size $\dz$ and (ii) the length
$T$ of the time window (recall that time $t$ in \eqref{e1_01} is a spatial coordinate),
and also (iii) a violation of the so called principle of ``frozen coefficients"
(see below).  Property (i) has been
considered both numerically and analytically for the monochromatic wave background, and
numerically for the soliton background, in \cite{J}. We will explain later on that this
property, i.e. the high sensitivity of the instability to the step size $\dz$, 
has different dependence on the time window length $T$ and also requires a different
mathematical description, for the monochromatic and soliton backgrounds.
As for properties (ii) and (iii), they, to the best of our knowledge, 
have not been previously systematically studied for any numerical scheme.

In Section 2, we will present results of our numerical simulations of Eq.~\eqref{e1_01} 
by the split-step Fourier method \eqref{e1_02} on the background of the soliton
\be
u_{\rm sol}=A\,\sqrt{\frac2{\gamma}}\,\sech\Big( \frac{A\,t}{\sqrt{-\beta}} \Big)\,
\exp\big(iK z\big)\, \equiv \, U(t)e^{iKz}, \qquad K=A^2\,.
\label{e1_08}
\ee
These will illustrate the unusual instability properties listed in the previous
paragraph. In Section 3, we will develop an analytical theory of the instability of the
slit-step Fourier method on the background of the soliton which explains all these observed
propeties. A comparison of this theory with numerical simulations
is presented in Section 4.
In Section 5 we summarize this work. Appendices 1 and 2 contain auxiliary results, and
Appendix 3 recovers the results of Weideman and Herbst \cite{WH}, illustrated in 
Fig.~\ref{fig_1}, via the analytical method presented in Section 3.

% --------------------------------------------------------------------------------------------------------
\section{Numerical study of the instability of the split-step Fourier method on the background of a
soliton }
\setcounter{equation}{0}

We numerically simulated Eq.~\eqref{e1_01} with $\beta=-1$, $\gamma=2$ on the background of the
soliton \eqref{e1_08} with $A=1$ using algorithm \eqref{e1_02}. We added to the initial soliton
a very small white (in $t$) noise
component, which served to reveal unstable Fourier modes sooner than if they had developed from
the round-off error. Thus, the initial condition for our numerical experiments was
\be
u_0(t)=A\,\sech (At) + \xi(t), \qquad A=1,
\label{e2_01}
\ee
and $\xi(t)$ was a Gaussian random process with zero mean and the standard deviation $10^{-10}$.

\begin{figure}[h]
\vspace{-1.6cm}
\rotatebox{0}{\resizebox{7cm}{9cm}{\includegraphics[0in,0.5in]
 [8in,10.5in]{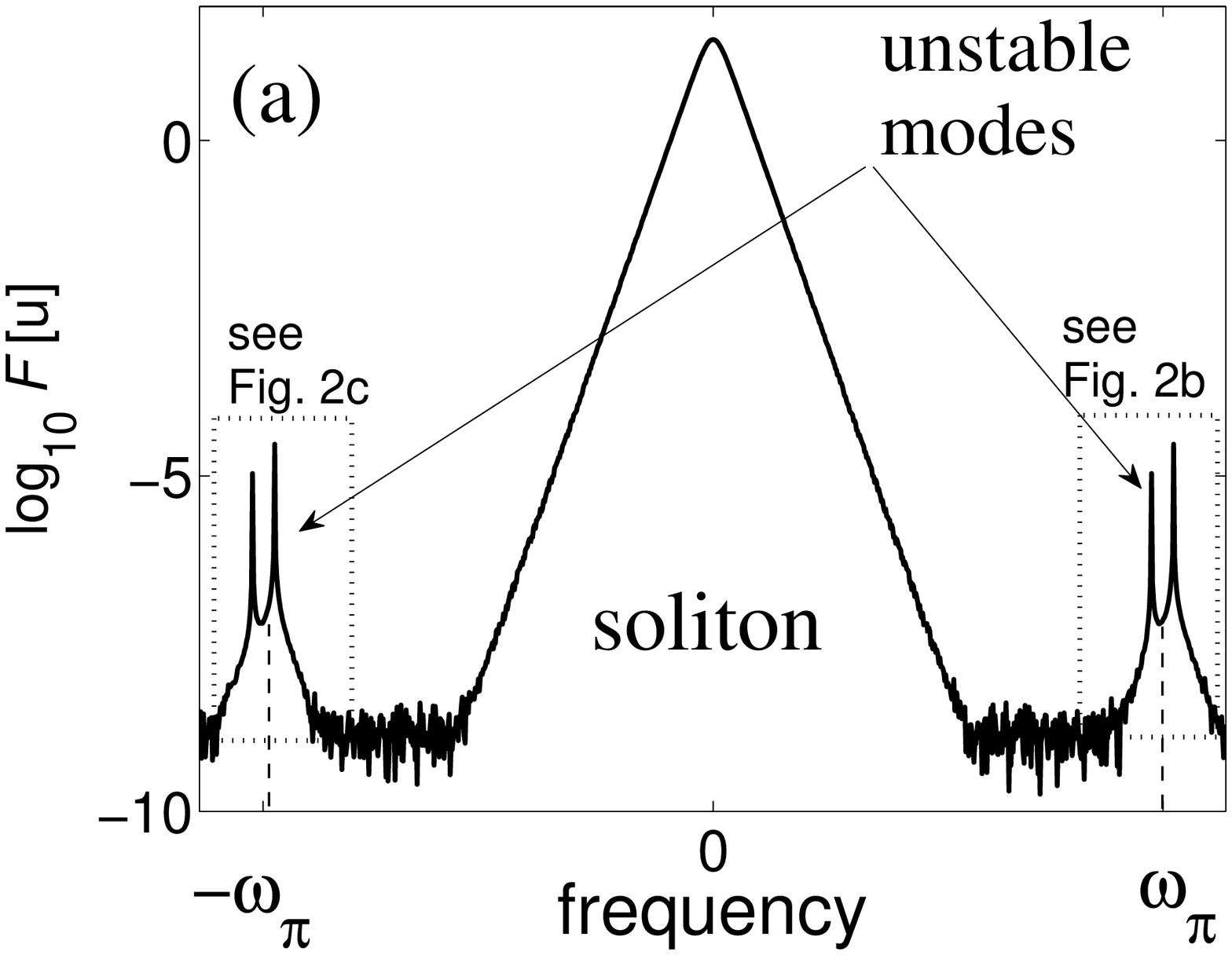}}}

\vspace{-2.5cm}

\mbox{ 
\begin{minipage}{7cm}
\rotatebox{0}{\resizebox{7cm}{9cm}{\includegraphics[0in,0.5in]
 [8in,10.5in]{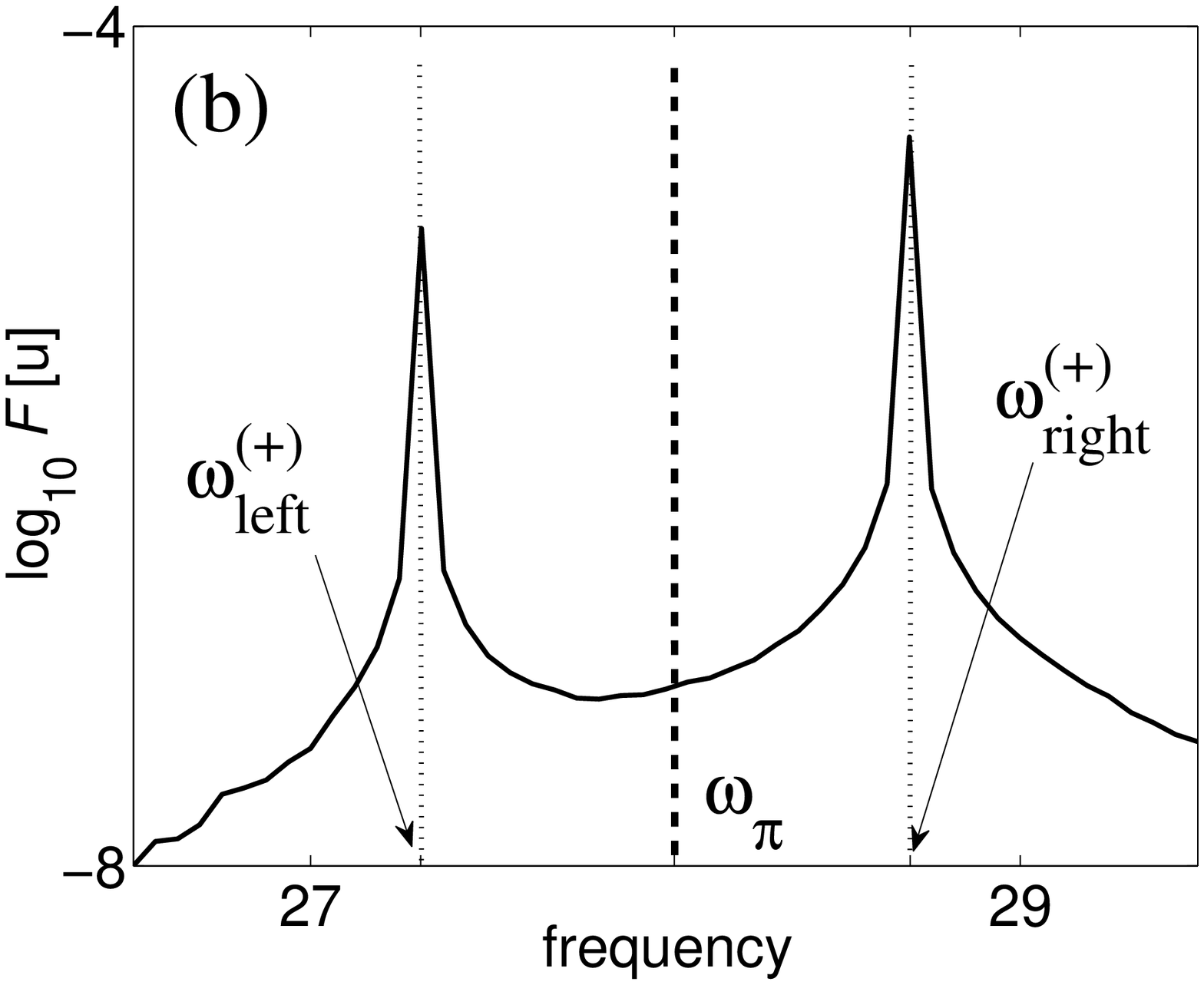}}}
\end{minipage}
\hspace{0.1cm}
\begin{minipage}{7cm}
\rotatebox{0}{\resizebox{7cm}{9cm}{\includegraphics[0in,0.5in]
 [8in,10.5in]{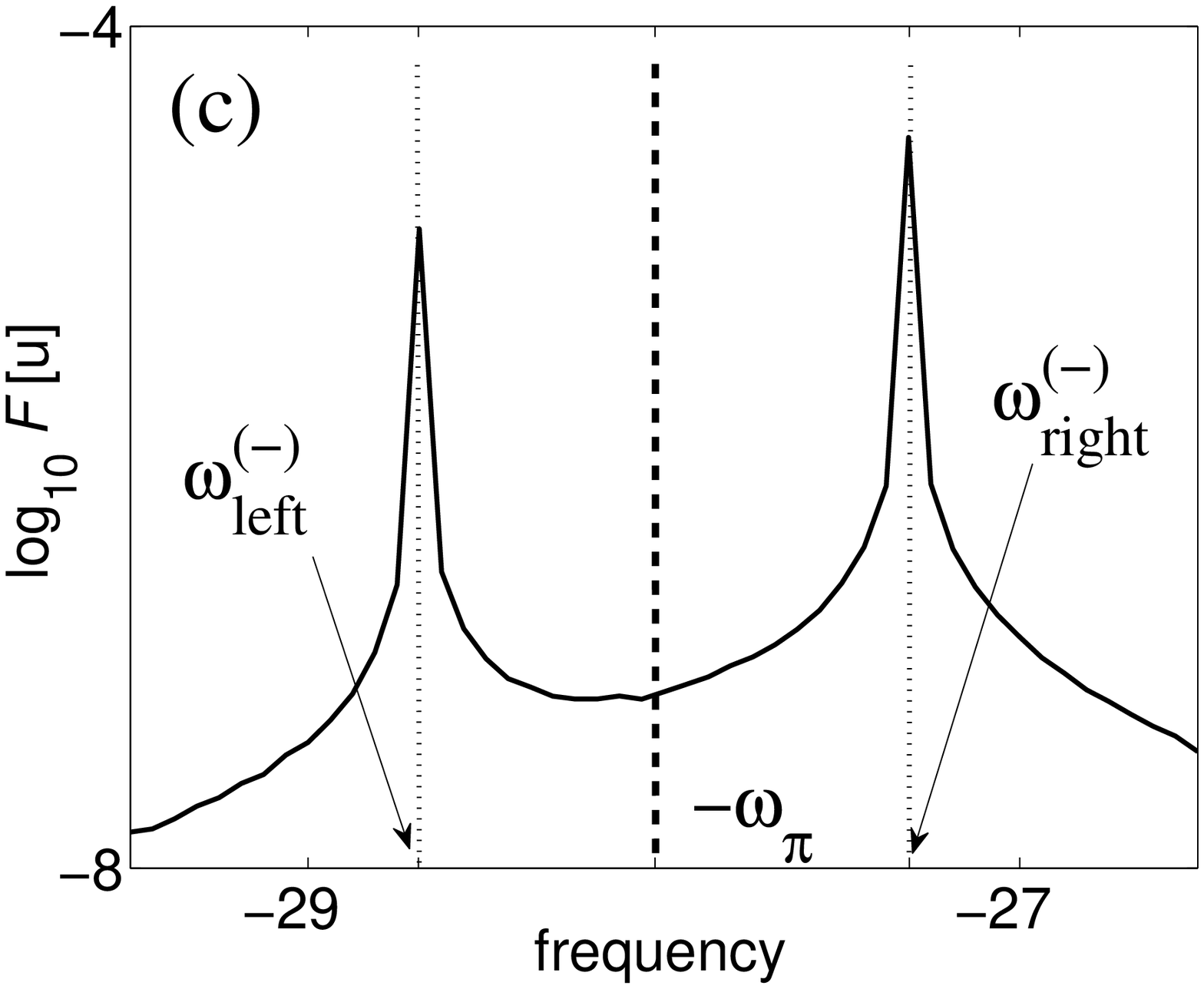}}}
\end{minipage}
 }
\vspace{-1.6cm}
\caption{The spectrum of the solution of \eqref{e1_01} and \eqref{e2_01} at $z=500$ with $\dz=0.0040$,
computed by method \eqref{e1_02}. Other parameters are listed in the text after Eq.~\eqref{e2_01}. 
}
\label{fig_2}
\end{figure}

To begin, we considered the spatial grid \ $-T/2 < t \le T/2$ with $2^{10}$ grid points and width
$T=32\pi$, i.e. about two orders of magnitude wider than the soliton. This results in
the spectral grid being the interval \ $-32 < \omega \le 32$ with the frequency spacing of
$\Delta\omega=2\pi/T=0.0625$. 
As expected, we did not observe any instability for $\dz < \dz_{\rm thresh}\approx 0.0031$.
Above the threshold, we ran the simulations with $\dz$ ranging from $0.004$ to $0.006$ with the increment
of $0.0001$ (i.e., $\dz=0.0040, \,0.0041, \, \ldots \, 0.0060$) up to the maximum distance of
$z_{\max}=500$ and observed instability only at a few values of $\dz$ listed in Table 1.
A typical Fourier spectrum of the unstable solution at $z_{\max}$ is illustrated in Fig.~\ref{fig_2}
for $\dz=0.0040$. 
The instability increment (see \eqref{e1_05b}) listed in Table 1 was computed as
\be
{\rm Re}(\lambda) = \frac{ 
     \left( \mbox{peaks' exponent} \right) \,-\, \left( \mbox{noise's exponent} \right) }{z_{\max}}
     \,\ln 10\,,
\label{e2_02}
\ee
where the peaks's exponent refers to the average of the decimal
exponents of the two peaks in Fig.~\ref{fig_2}(b),
and the exponent corresponding to the average noise level was estimated to be $-8.8$ for this 
set of simulations. Estimate \eqref{e2_02} may not be very accurate, as the so computed increment depends
on the noise levels at the frequencies of the unstable Fourier modes (hence the peaks in Fig.~\ref{fig_2}
are slightly different), but it still provides a reasonable measure of the instability rate.

\begin{table}
\bc % \footnotesize
\begin{tabular}{|c|c|c|c|c|}  \hline 
 $\dz$   &  $\omega^{(+)}_{\rm left}-\opi$  &  $\opi$  &  $\omega^{(+)}_{\rm right}-\opi$ 
  & Re$(\lambda)$  \\ \hline 
 0.0040  &  -0.72  &  28.03  &  0.66  &  0.019  \\ \hline 
 0.0048  &  -0.39  &  25.58  &  0.35  &  0.031  \\ \hline
 0.0054  &  -0.62  &  24.12  &  0.57  &  0.022  \\ \hline
 0.0055  &  -0.52  &  23.90  &  0.48  &  0.024  \\ \hline
 0.0058  &  -0.46  &  23.27  &  0.42  &  0.022  \\ \hline                                                  
\end{tabular}
\ec
\caption{Parameters of the unstable frequencies' peaks when $T=32\pi$, number of grid points is
$2^{10}$, and $\dz$ is varied as \ $0.0040,\,0.0041,\,0.0042,\, \ldots\,, 0.0060$.
The notations $\omega^{(+)}_{\rm right,\;left}$ are introduced in Fig.~\ref{fig_2}.
 }
\end{table}

From Fig.~\ref{fig_2} we first observe that the instability peaks at negative frequencies
are not reflectionally symmetric relative to such peaks at positive frequencies, in contrast to the
instability peaks on the background of a monochromatic wave (see the caption to Fig.~\ref{fig_1}).
Rather, for all the simulations we ran,
the negative-frequency peaks appear to be a {\em shifted} replica of the positive-frequency peaks.
Moreover, the frequencies of the left and right peaks on the same side of $\omega=0$ are 
slightly asymmetric with respect to $\opi$: see Table 1. To further investigate how the
peaks are related to one another, we ran a simulation where we placed a narrow filter that totally
suppressed the field around $\omega=\omega^{(+)}_{\rm right}$. As a result, the peak at
$\omega=\omega^{(-)}_{\rm right}$ was suppressed but the peaks at $\omega=\omega^{(\pm)}_{\rm left}$
remained intact. This corraborated the aforementioned observation that the pairs of peaks
appear to be shifted replicas of each other rather than reflected about $\omega=0$. 
In light of this, it might be somewhat surprising that the peak's frequencies are found to
be related by a reflectional symmetry:
\be
\omega^{(-)}_{\rm left} = - \omega^{(+)}_{\rm right}, \qquad
\omega^{(-)}_{\rm right} = - \omega^{(+)}_{\rm left}.
\label{e2_03}
\ee
All these observations are explained by the theory in Section 3.

Another conspicuous observation, made from Table 1, is that there seems to be no apparent order
in the frequencies and heights of the instability peaks as functions of $\dz$. To expand
on this, we zoomed in on the interval $0.0045 \le \dz \le 0.0046$, at the end points of which
we had found {\em no} instability. We varied $\dz=0.00451,\,0.00452,\, \ldots \,\, 0.00459$ and
found that the instability occured at $\dz=0.00451$ and persisted up to $\dz=0.00456$, with its
increment Re$(\lambda)$ gradually decreasing from $0.036$ to $0.014$ and the inter-peak
spacing $\omega^{(+)}_{\rm right}-\omega^{(+)}_{\rm left}$ gradually increasing from $0.25$ to
$1.56$. A weaker instability, with the inter-peak spacing of about $2.0$, was also observed at
$\dz=0.00459$; note that no instability was observed at $\dz=0.00457$ and $0.00458$.

One could argue that such an irregular dependence of the instability on $\dz$ occurs simply
because the width of the band of unstable modes is just slightly greater than the mesh size $\Delta\omega$
of the frequency grid.
Indeed, the unstable band's width estimated from the monochromatic-background
case is $\gamma A^2 /(|\beta|\opi) \approx 0.07$ for the parameters listed above
(see the caption to Fig.~\ref{fig_1} and formula \eqref{a3_06} in Appendix 3),
while $\Delta\omega\equiv 2\pi/T = 0.0625$. 
In such a case, the instability features should be strongly sensitive to where
inside the instability band the frequency $2\pi\ell/T$ falls. This was pointed out in \cite{J};
see also the end of Appendix 3 below.
Then one would expect that the aforementioned sensitivity is to be alleviated by taking a wider time window,
because then the frequency mesh size $\Delta\omega$ would decrease and more frequencies from the
numerical grid would fall into the unstable band.

\begin{figure}[h]
\vspace*{-1cm}
\mbox{ 
\begin{minipage}{7cm}
\rotatebox{0}{\resizebox{7cm}{9cm}{\includegraphics[0in,0.5in]
 [8in,10.5in]{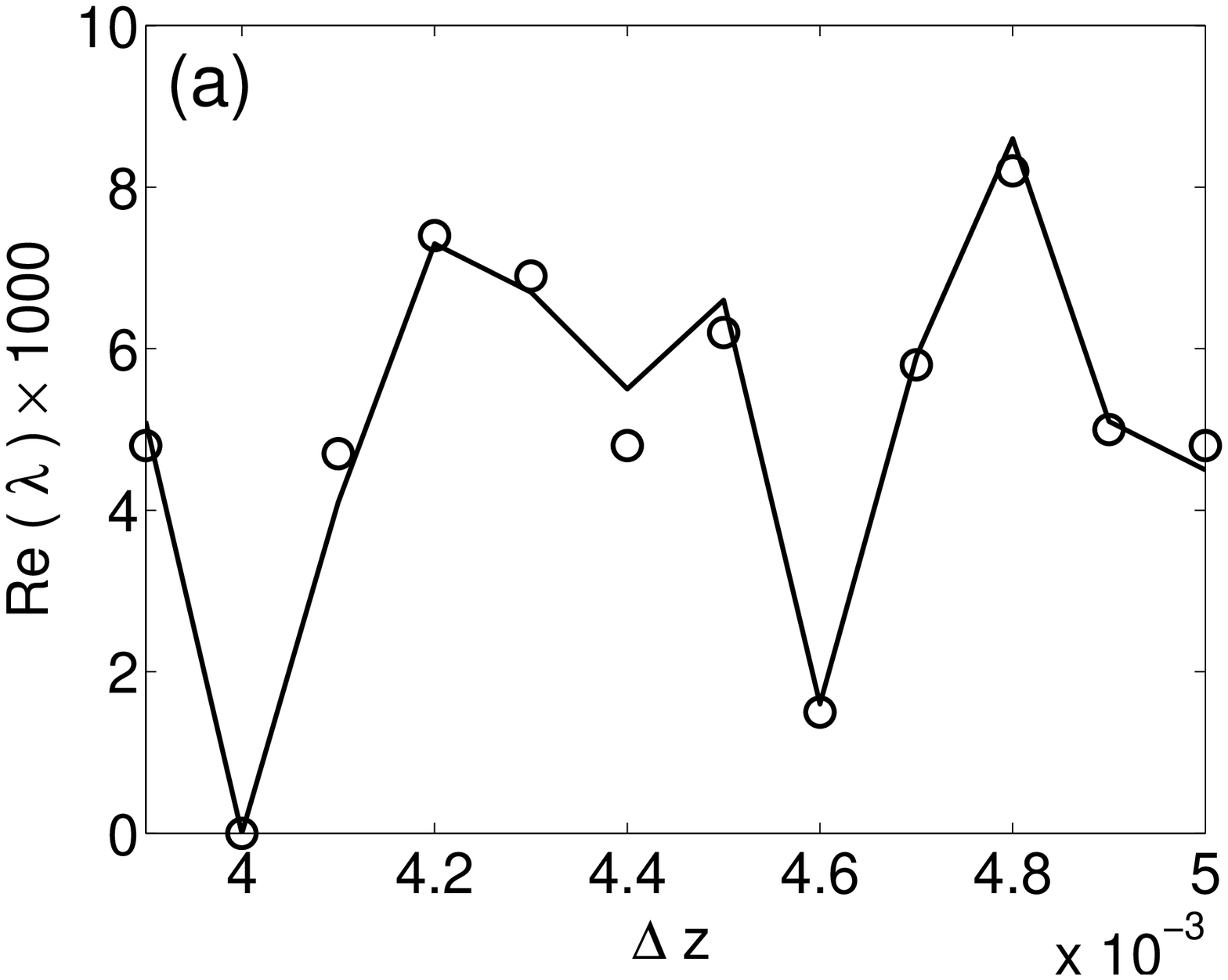}}}
\end{minipage}
\hspace{0.1cm}
\begin{minipage}{7cm}
\rotatebox{0}{\resizebox{7cm}{9cm}{\includegraphics[0in,0.5in]
 [8in,10.5in]{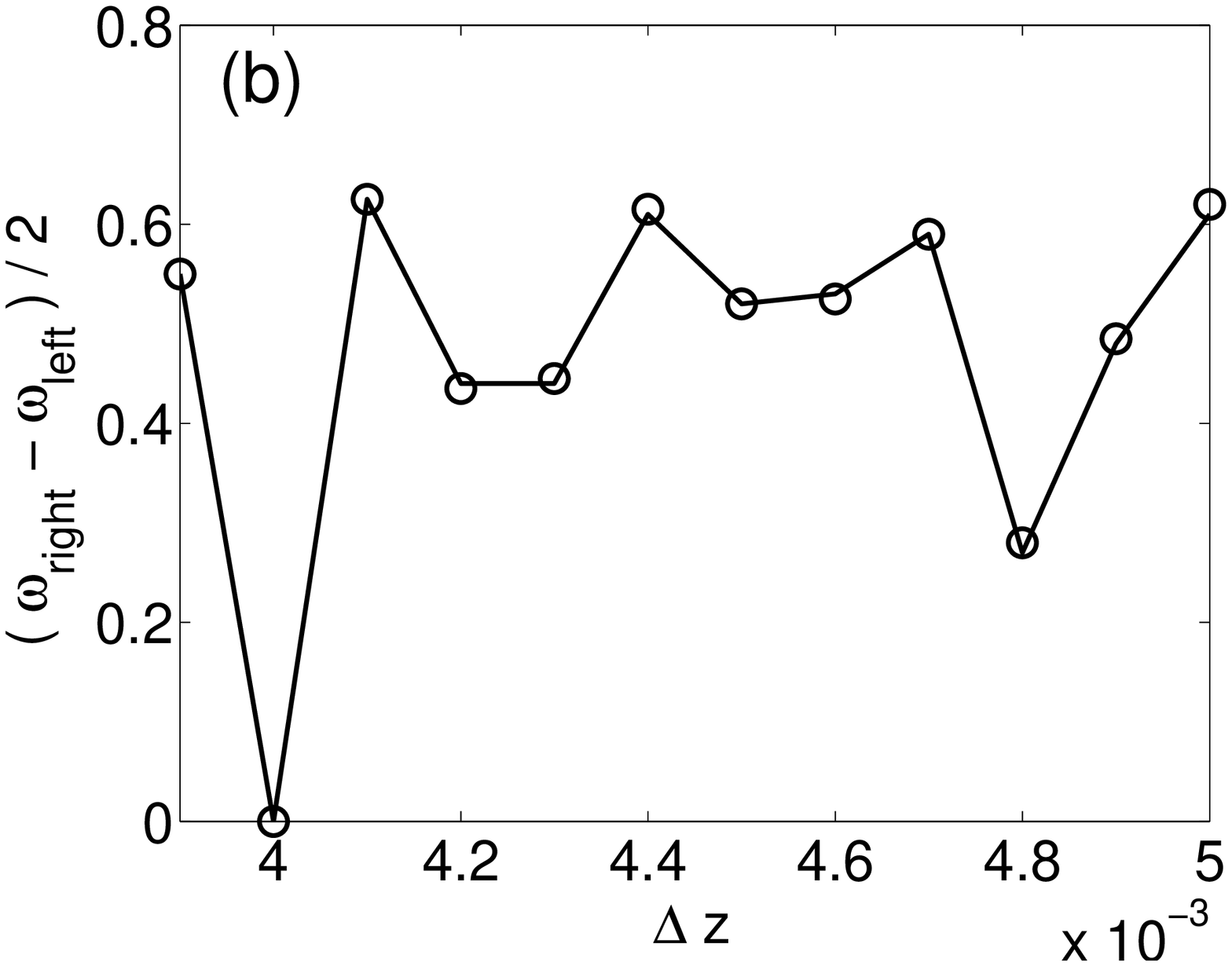}}}
\end{minipage}
 }
\vspace{-1.6cm}
\caption{Instability increment and unstable peaks' half-separation at $z=2000$,
as a function of the step size $\dz$. The initial condition is \eqref{e2_01}.
The numerical domain has length $T=128\pi$ and contains $2^{12}$ grid points.
Note that for $\dz=0.0040$, where no instability is found, we {\em defined}
the corresponding $\omega_{\rm right}-\omega_{\rm left}=0$. 
The open circles and solid line correspond, respectively, to the numerical results 
and to the results of analytical calculations reported in Section 4.2.
}
\label{fig_3}
\end{figure}

\begin{figure}[h]
\vspace*{-1cm}
\mbox{ 
\begin{minipage}{7cm}
\rotatebox{0}{\resizebox{7cm}{9cm}{\includegraphics[0in,0.5in]
 [8in,10.5in]{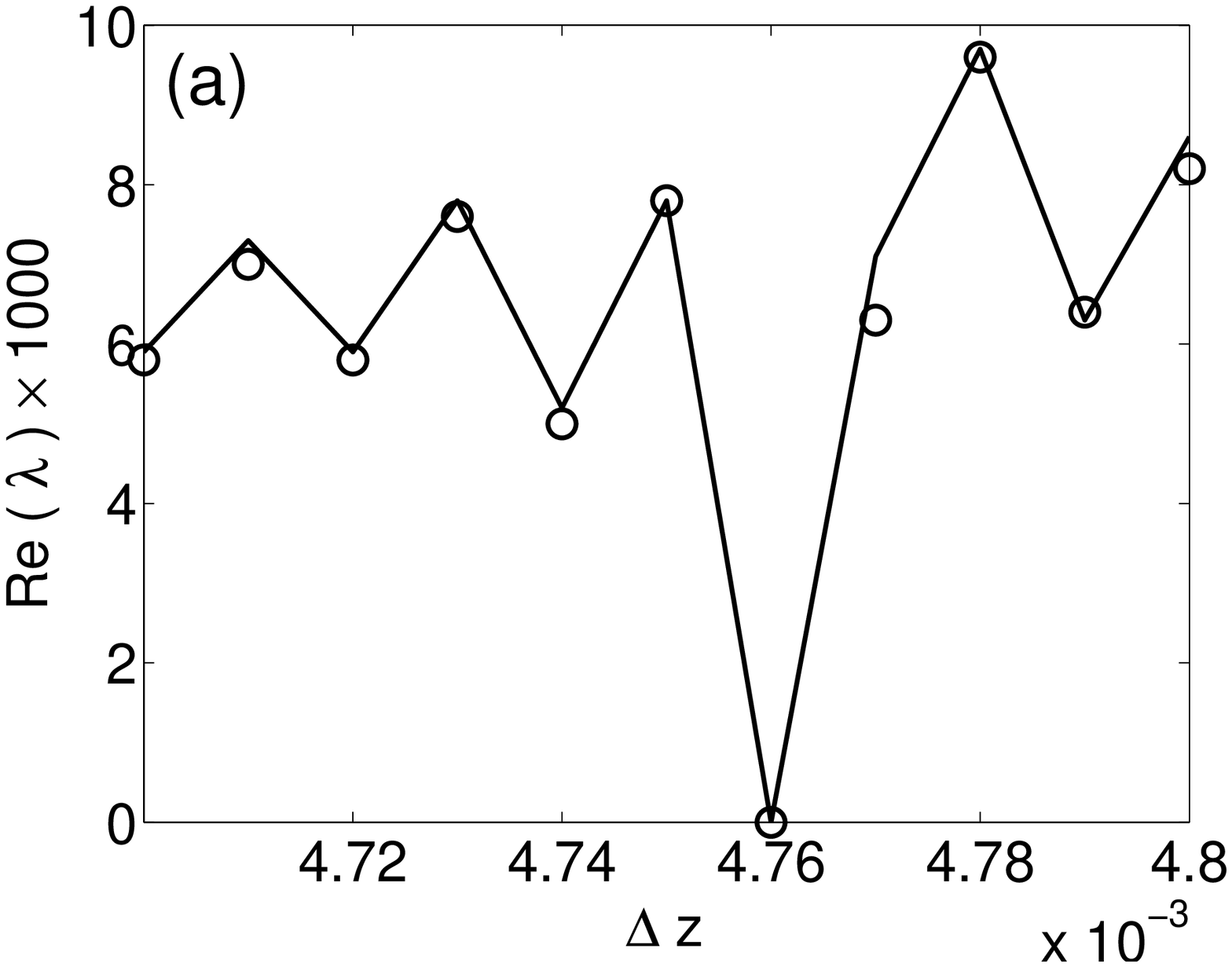}}}
\end{minipage}
\hspace{0.1cm}
\begin{minipage}{7cm}
\rotatebox{0}{\resizebox{7cm}{9cm}{\includegraphics[0in,0.5in]
 [8in,10.5in]{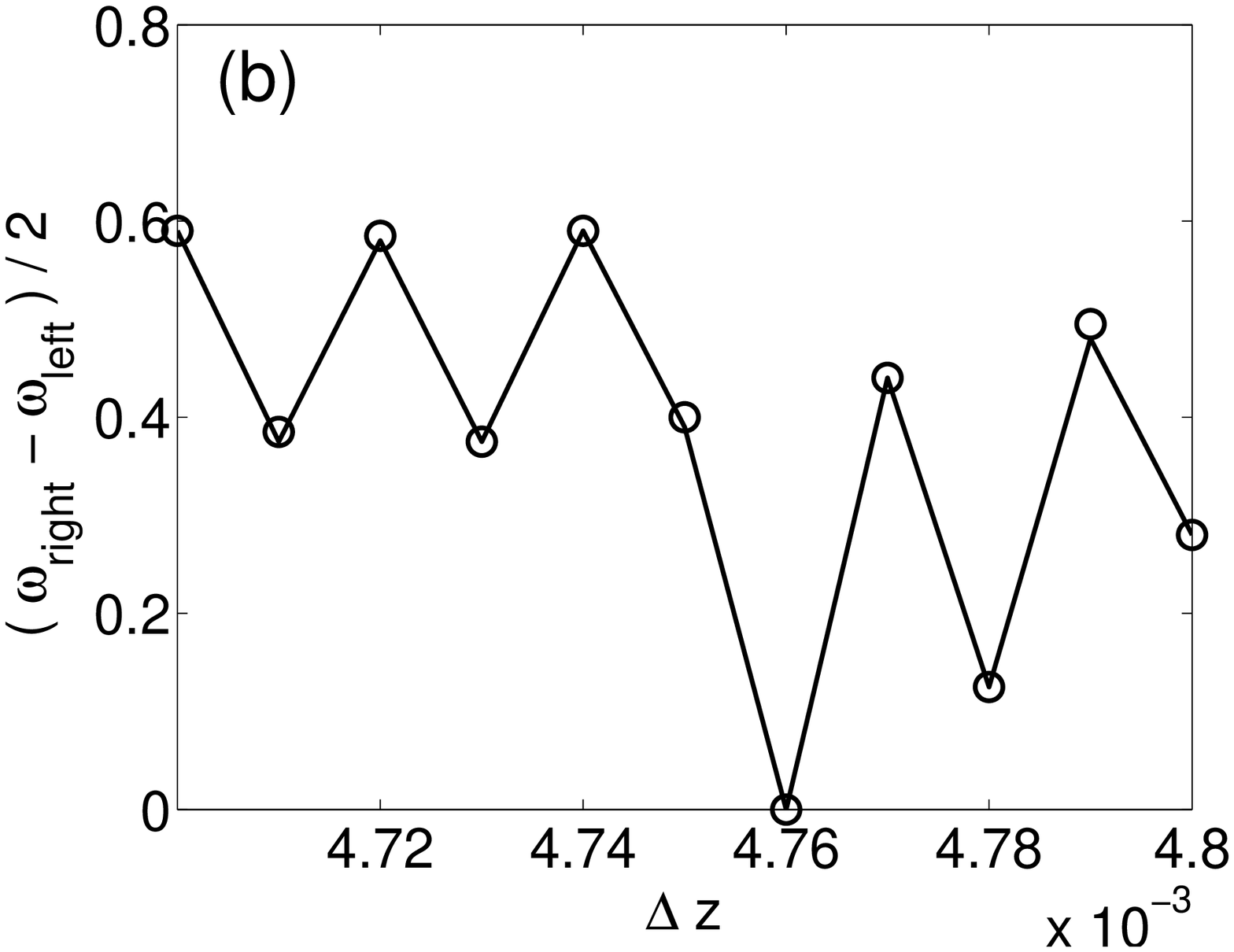}}}
\end{minipage}
 }
\vspace{-1.6cm}
\caption{Same as in Fig.~\ref{fig_3}, but for a different interval of $\dz$ values.
}
\label{fig_4}
\end{figure}

\begin{figure}[h]
\vspace*{-1cm}
\mbox{ 
\begin{minipage}{7cm}
\rotatebox{0}{\resizebox{7cm}{9cm}{\includegraphics[0in,0.5in]
 [8in,10.5in]{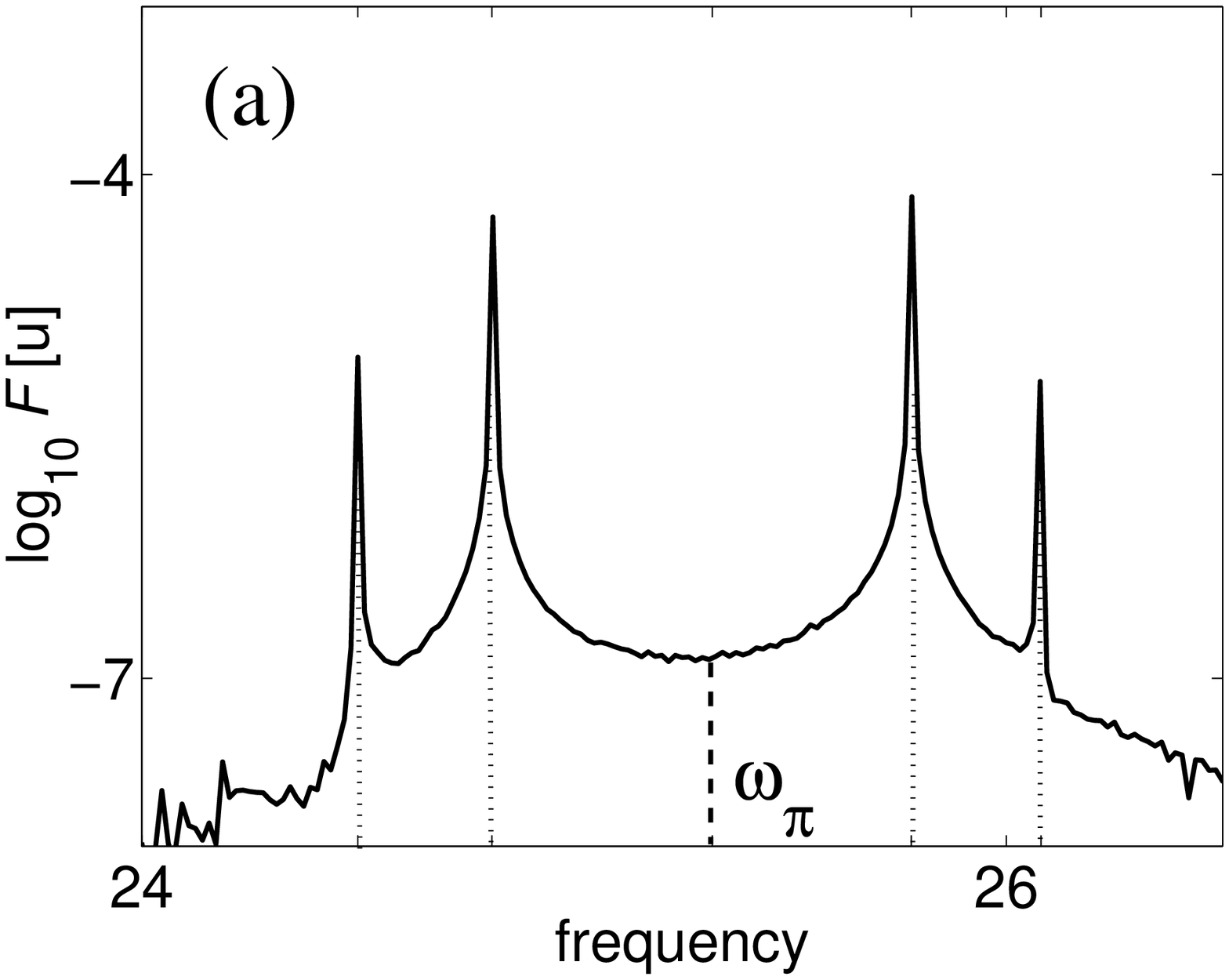}}}
\end{minipage}
\hspace{0.1cm}
\begin{minipage}{7cm}
\rotatebox{0}{\resizebox{7cm}{9cm}{\includegraphics[0in,0.5in]
 [8in,10.5in]{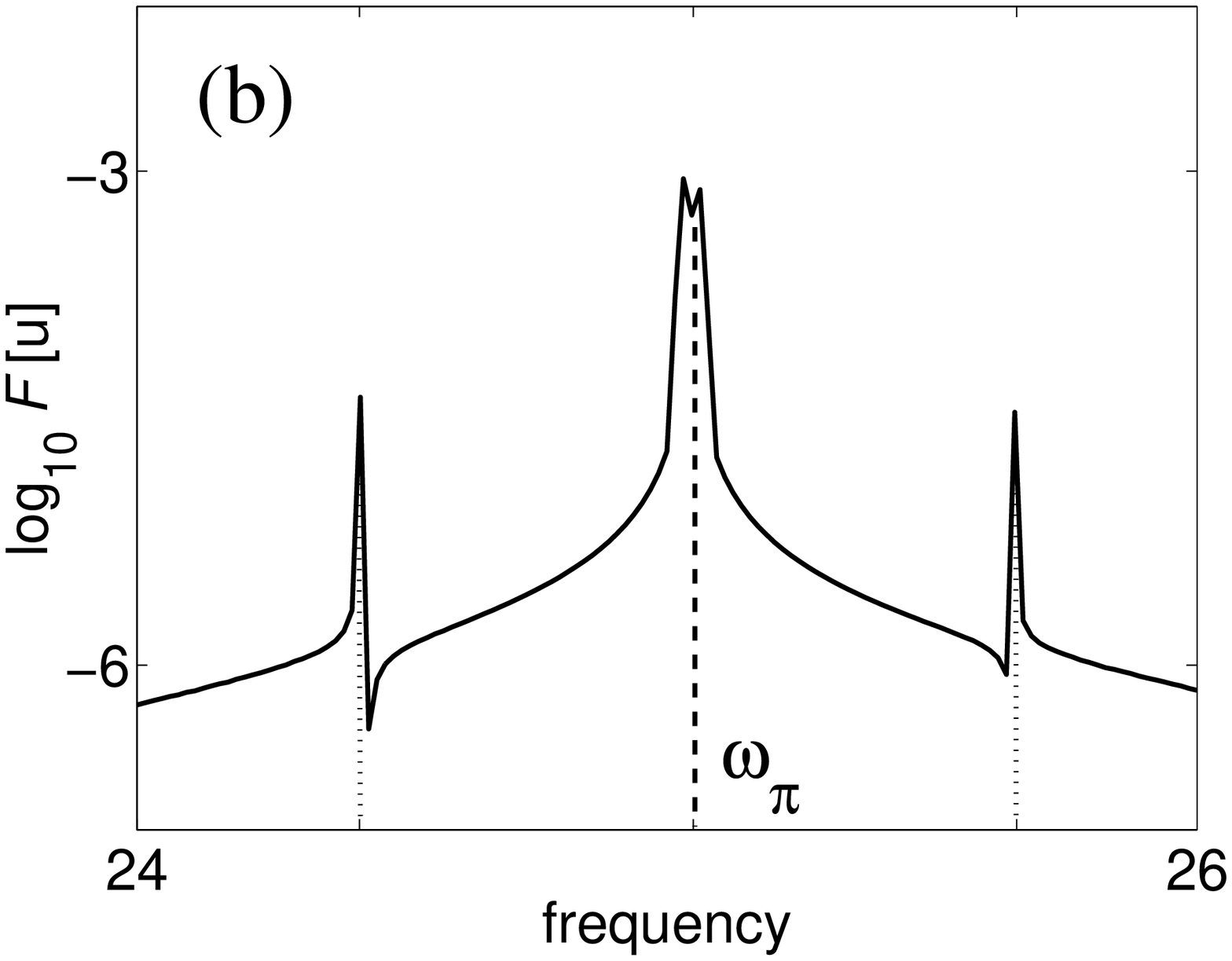}}}
\end{minipage}
 }
\vspace{-1.6cm}
\caption{``Less typical" instability spectra. Both panels are for
initial condition \eqref{e2_01}, the propagation distance $z=2000$,
and the numerical domain with $T=128\pi$ and $2^{12}$ grid points.
 \ (a) \ $\dz=0.0049$; similar spectra are found for $\dz=0.0044, \,
0.0045, \, 0.00474$; \ (b) \ $\dz=0.0050$.
}
\label{fig_5}
\end{figure}

To verify the validity of such an argument, we re-run the
simulations described above taking a four times wider time window, i.e. $T=128\pi$,
while retaining the same $\Delta t$ (thus quadrupling the number of grid points). For these
new simulations, we estimate that instability peaks should contain about four grid points
and hence may expect that the high sensitivity to the value of $\dz$ reported above is to be alleviated.
What we find instead is {\em an opposite} of this statement: the dependence of the locations and heights
of the instability peaks remains at least as irregular as for the smaller value of $T$, but now
the instability is observed ``more often" than in Table 1. In Figs.~\ref{fig_3} and \ref{fig_4} 
we plot the observed values of the instability increment and half of the inter-peak spacing,
$(\omega^{(+)}_{\rm right}-\omega^{(+)}_{\rm left})/2$, when $\dz$ is varied between $0.0039$
and $0.0050$ with step $0.0001$ and between $0.00471$ and $0.00479$ with step $0.00001$.
Since the instability rates have now been found to be
about a factor of four lower than in Table 1, we had to
run our simulations up to a greater distance, $z_{\max}=2000$. Note that the
instability characteristics reported for $0.0040 \le \dz \le 0.0050$ in Table 1 do not
match those shown in Fig.~\ref{fig_3}.
Moreover, we find that the spectra of unstable modes may look
qualitatively different than in Fig.~\ref{fig_2}. Namely, these spectra for $\dz=0.0044, \,
0.0045, \, 0.00474, \, 0.0049$, look like the one shown in Fig.~\ref{fig_5}(a),  while for
$\dz=0.0050$ it is shown in Fig.~\ref{fig_5}(b). Also, contrary to our expectation,
we observe that in most cases the instability peaks still contain only one grid point;
exceptions are the central peak for $\dz=0.0050$ and the peaks for $\dz=0.00478$, which are
spaced very closely. This fact will be emphasized when we describe a challenge in the
analytical description of the instability in Section 3.

We have already mentioned that keeping the spatial mesh size $\Delta t$ intact but increasing
the spatial window's length $T$ by a factor of four considerably changes parameters of the instability. 
We now show that just {\em slightly} changing
$T$ (and hence correspondingly slightly changing $\Delta t$) may change instability
parameters {\em dramatically}. 
Such a sensitivity to the length of the spatial window is
not observed for finite-difference methods.
Note that $T=128\pi \approx 402.1$. In Table 2 we list the instability parameters
observed when we decrease $T$ by about 1\% or even less. Again, these results appear to
be totally irregular. Similarly to the situation mentioned at the end of the previous paragraph, 
only the closely spaced instability peaks at $T=398$ contain several grid points; the peaks for
the other values of $T$ each contain only one grid point.

\begin{table}
\bc  % \footnotesize
\begin{tabular}{|c|c|c|c|}  \hline 
 $T$   &  $\omega^{(+)}_{\rm left}-\opi$   &  $\omega^{(+)}_{\rm right}-\opi$ 
  & Re$(\lambda)$  \\ \hline 
 397  &  -0.44  &  0.40  &  0.0062  \\ \hline 
 398  &  -0.10  &  0.08  &  0.0096  \\ \hline
 399  &  -1.09  &  1.02  &  0.0024  \\ \hline
 400  &  -0.33  &  0.29  &  0.0076  \\ \hline
 $128\pi$  &  -0.39  &  0.36  &  0.0076  \\ \hline                                                  
\end{tabular}
\ec
\caption{Parameters of the unstable frequencies' peaks when $\dz=0.0043$ ($\opi=27.03$), 
number of grid points is $2^{12}$, and $T$ is varied.
 }
\end{table}

Finally, let us show that the so called principle of ``frozen coefficients" 
\cite{VonNeumannRichtmeyer}, which is
known to apply to finite-difference schemes, 
does {\em not} hold for the split-step Fourier method \eqref{e1_02} on the 
background of a soliton. This principle says the following. Suppose one has an
evolution equation with a spatially varying coefficient, say $c(t)$
(recall that in this paper, $t$ is a spatial variable, 
while $z$ is the evolution variable). Near each $t=t_0$,
this coefficient can be approximated by $c(t_0)$. Then one can apply the standard von Neumann
stability analysis to the equation with the {\em constant} (``frozen") coefficient $c(t_0)$. 
Then the scheme is deemed unstable if such an analysis reveals instability for at least
one value of $c(t_0)$.
This principle works quite well for finite-difference schemes.
An intuitive explanation fot this is that the unstable modes usually have high spatial frequency,
and so they ``see" a relatively slowly-varying coefficient $c(t)$ as being approximately constant
near each $t=t_0$. For the split-step Fourier method, the unstable modes also have high
spatial frequencies, $\approx \pm\opi$. Yet, the principle of ``frozen coefficients" 
does not apply to this method, as we illustrate below.

Indeed, we have already seen that the instability of method \eqref{e1_02} on the background
of a soliton is different from that on the background of a monochromatic wave of the same
amplitude.
This fact was originally stated by Weideman and Herbst (see the last section in \cite{WH}).
Now we will present three examples that deal with a multi-soliton background.
Conclusions from these examples are at odds with our intuition based on the experience
with finite-difference schemes, and they reveal yet another way in which the principle
of ``frozen coefficients" is violated for method \eqref{e1_02}.
As at the beginning of this Section, let us 
use the numerical domain with $T=32\pi \approx 100.5$ and $2^{10}$ grid points, 
but instead of the {\em single}-soliton initial condition \eqref{e2_01}, consider three
{\em well-separated} solitons:
\be
u_0(t)=\sech (t+33.5) + \sech (t) + \sech (t-33.5) \,+\, \xi(t),
\label{e2_04}
\ee
where $\xi(t)$ is the same as in \eqref{e2_01}. The spacing between the solitons
is chosen to be sufficiently large so as to avoid their interaction. According to the
principle of ``frozen coefficients", the instability with this initial condition must be 
the same as with \eqref{e2_01}. However, it turns out to be quite different.
First, for $\dz=0.0040$, when there was an instability (see Table 1) for the initial condition
\eqref{e2_01}, there is {\em no} instability for initial condition \eqref{e2_04}.
Second, for $\dz=0.0052$, when there was {\em no} instability according to Table 1,
now there is a strong instability with the increment of about $0.10$ (i.e., more than three
times greater than the largest increment in Table 1). Third, for $\dz=0.0058$, the 
instability was observed with both initial conditions \eqref{e2_01} and \eqref{e2_04},
but for the latter case it was more than three times as strong (having the increment of $0.072$).

\smallskip

To summarize, in this Section we have presented results of numerical simulations of the
split-step Fourier method \eqref{e1_02} with soliton initial conditions. These results
illustrate the following features that are drastically distinct from features exhibited
by unstable finite-difference schemes:
\begin{itemize}
\item
The spectra of the unstable modes exhibit strong and irregular sensitivity to the
values of the step size $\dz$ and the length of the spatial window $T$;
\item
Instability on the background of several well-separated --- and hence non-interacting ---
pulses can be drastically different from the instability on the background of a single
pulse (or a different number of pulses). 
\end{itemize}
The theory that we will present in the next section is able to {\em quantitatively}
explain all these features.

% -----------------------------------------------------------------------------------------------
\section{Analytical theory of the instability of the split-step Fourier method on the background of a
soliton }
\setcounter{equation}{0}

\subsection{Key idea and main challenge}

First, we will show that a straightforward application of the
von Neumann analysis is unlikely to provide analytical insight about stability or
instability of the split-step method on the soliton background. Then we will outline
the idea of an alternative approach and point out the mathematical challenge that such an
aproach would have to resolve.

Substituting expression \eqref{e1_05a} with $\bar{u}=\usol$ into
algorithm \eqref{e1_02} and keeping only terms linear in $\tu_n$, one obtains after
taking the Fourier transform:
\be
\F[\tu_{n+1}] = e^{i\beta\omega^2 \dz} \, 
\F \left[ e^{i\gamma |\usol|^2 \dz} \big( \tu_n + i\gamma\dz ( \usol^2\tu_n^* + |\usol|^2 \tu_n) \big) 
 \right]\,,
\label{e3_01}
\ee
where $\usol$ is given by \eqref{e1_08}. The right-hand side of \eqref{e3_01} 
describes coupling of Fourier modes \ $\F[\tu_n](\omega_{\ell})$ \ with different
$\omega_{\ell}=2\pi\ell/T$, $\ell=0,\,\pm 1,\,\pm 2, \, \ldots $ via, e.g., the
convolution term \ $\F\big[ |\usol|^2\tu_n \big]$, because $|\usol|$ is a function of the spatial variable $t$.
(Note that this problem did not occur in the stability analysis on the background of a monochromatic
wave \eqref{e1_04} since there $|\bar{u}|=|u_{\rm cw}|$ is $t$-independent 
and hence the corresponding
equation, studied by Weideman and Herbst, coupled only the modes $\omega_{\ell}$
and $-\omega_{\ell}$.) Solving for eigenvalues of such a coupled multi-mode
system, while numerically
feasible, would unlikely provide any insight of how the instability can occur.
Such an insight is provided by an alternative analytical approach described below.

The key idea of this approach comes from the instability spectra shown in
Figs.~\ref{fig_1}, \ref{fig_2}, and \ref{fig_5}. Namely, we note that only narrow bands of Fourier
modes with frequencies near $\pm\opi$ can become unstable. The reason for this
will be presented later. Therefore, it is appropriate
to seek the numerical error $\tu_n$, defined in \eqref{e1_05a}, as consisting of two 
quasi-monochromatic waves whose carrier frequencies are approximately $\pm\opi$. 
These waves can become unstable via their interaction mediated by their scattering on the soliton.
Since these are high-frequency waves (see \eqref{e1_06} and note that $\dz\ll 1$), they ``see"
the soliton, whose temporal width is $O(1)$, as a narrow, and hence small, perturbation. 
This is one of the reasons
that the instability is weak, as seen in Section 2.
Using the weakness of the instability, we will approximate its evolution by a differential,
rather than difference, equation. Analysis of the former kind of equation
is considerably easier than that of the latter one.

%  This idea may qualitatively explain the high sensitivity of the instability to the length
%  $T$ of the spatial window, the value of $\opi$ (equivalently, to $\dz$), and the number and
%  arrangement of scatterers (i.e., the solitons), as summarized at the end of Section 2.
%  Indeed, all these features affect the 

The main mathematical challenge that this approach needs to address is a {\em three}-scale
nature of this problem, as can be seen from Fig.~\ref{fig_2}. One scale is set by $\opi \gg 1$,
where the strong inequality is a consequence of $\dz \ll 1$ (see \eqref{e1_06}), 
as is usually the case in
simulations. Accordingly, let us define a small parameter
\be
\epsilon=\frac1{\opi} \ll 1, \qquad \epsilon=O\big(\sqrt{\dz}\big),
\label{e3_02}
\ee
which will play a prominent role in what follows. In addition to this $O(\epsilon^{-1})$ scale,
there is also a scale $O(1)$: The separation between $\opi$ and the frequencies of the most unstable
modes, seen as peaks in Fig.~\ref{fig_2}, appears to be of this order of magnitude. Finally, the third scale 
is set by the spectral width of the instability peaks. Indeed, as noted in Section 2, even for the highest
spectral resolution reported there, those peaks, in most cases, contain only one grid point.
(The pedestals seen around these peaks can be shown to occur due to modulational instability on
the background of the peaks, i.e. are not directly caused by a numerical instability.) 
Based on the information provided by Figs.~\ref{fig_2} and \ref{fig_5},
it is even impossible to tell whether this third scale is determined by $\epsilon$. We will show below
that it is rather determined by several parameters of the problem.
%  and, most prominently, by the length of the temporal window $T$.
Incidentally, let us note that
the instability on the background of a monochromatic wave has only two scales: the location
$\pm\opi=O(\epsilon^{-1})$ of the instability peaks and the peaks' width, $O(\epsilon)$; see
Fig.~\ref{fig_1} and Appendix 3.

% ...............................................................................................

\subsection{Details of the theory}

First, let us simplify the right-hand side of \eqref{e3_01} by noting that in practice,
$\dz$ is always taken so as to guarantee $\gamma|\usol|^2\dz \ll 1$. Then, discarding
terms $O\big((\gamma|\usol|^2\dz)^2\big)$, one reduces \eqref{e3_01} to
\be
\F[\tu_{n+1}] = e^{i\beta\omega^2 \dz} \, 
\F \left[ \tu_n + i\gamma\dz ( \usol^2\tu_n^* + 2|\usol|^2 \tu_n) \right]\,.
\label{e3_03}
\ee
Let us now use the observation made in the previous subsection: As follows from Figs.
\ref{fig_1}, \ref{fig_2}, and \ref{fig_5}, the instability occurs only in narrow
spectral bands near $\pm\opi$, $\pm\omega_{2\pi}$, etc. We will focus on the instability
near $\pm\opi$; the analysis near $\pm\omega_{2\pi}$ etc. is similar. 
(In Appendix 2 we will show that
% if one carries out the subsequent analysis near an arbitrary frequency $\omega_0$, one still finds that
the instability can occur
{\em only} near $\pm\opi$, $\pm\omega_{2\pi}$, etc., thereby recovering the results presented
below.)
Accordingly,
%
%   let us expand the first term on the right-hand side of \eqref{e3_03} in powers of
%   \ $\beta(\omega^2-\opi^2)\dz=O(\epsilon)$, \ where the last order-of-magnitude estimate
%   follows from \eqref{e3_02} and the fact that we only consider frequencies satisfying
%   \ $|\omega-(\pm\opi)|=O(1)$; see the end od Section 3.1. The result of this expansion is:
%   %
%   \be
%   \hspace*{-7cm}
%   \F[\tu_{n+1}] = \, (-1) \, 
%   \F \left[ \tu_n + \left\{ i\beta(\omega^2-\opi^2)\dz\,\tu_n \right\} \, \right. 
%   \label{e3_04}
%   \ee
%   %
%   $$
%   \hspace*{2.5cm}
%   + \,
%   \left\{ i\gamma\dz ( \usol^2\tu_n^* + 2|\usol|^2 \tu_n) + 
%           \left(i\beta(\omega^2-\opi^2)\dz\right)^2 \tu_n  \right\} 
%   $$
%   $$
%   \hspace*{-2.1cm} \left. +  \, \left\{ \left(i\beta(\omega^2-\opi^2)\dz\right)^3 \tu_n  +  
%      \right.\right.
%   $$
%   $$
%   \hspace*{4cm} \left. 
%       i\beta(\omega^2-\opi^2)\dz\,\cdot\, i\gamma\dz ( \usol^2\tu_n^* + 2|\usol|^2 \tu_n) \Big\}
%       +  O(\epsilon^4) \right]\,.
%   $$
%
%   Here, the three consecutive sets of curly brackets group terms of orders $O(\epsilon)$, 
%   $O(\epsilon^2)$, and $O(\epsilon^3)$. The $(-1)$ in front of the right-hand side occurs due
%   to the definition of $\opi$, Eq.~\eqref{e1_06}. In Appendix 1 we show that if one carries out the
%   subsequent analysis with a similar expansion, but at a frequency that is not sufficiently close
%   to $\pm\opi$ (or $\pm\omega_{2\pi}$ etc.), then instability will not occur. 
%
let us rewrite \eqref{e3_03} as
\be
\F[\tu_{n+1}] = (-1)\cdot e^{i\beta(\omega^2-\opi^2) \dz} \, 
\F \left[ \tu_n + i\gamma\dz ( \usol^2\tu_n^* + 2|\usol|^2 \tu_n) \right]\,.
\label{e3_04}
\ee
where the $(-1)$ in front of the right-hand side occurs due
to the definition of $\opi$, Eq.~\eqref{e1_06}.
Note also that \ $\beta(\omega^2-\opi^2)\dz=O(\epsilon)$, \ which
follows from \eqref{e3_02} and the fact that we consider frequencies satisfying
 \ $|\omega-(\pm\opi)|=O(1)$; see the end of Section 3.1. Thus, \ 
 $\exp[i\beta(\omega^2-\opi^2) \dz]=1+O(\epsilon)\approx 1$  \ for the frequencies of interest.

To enable further analysis of the still intractable difference equation \eqref{e3_04} 
in the frequency domain,
let us convert it into a 
partial differential equation in the time domain. To this end, we first observe that for
the variable 
\be
\tv_n=(-1)^n \tu_n\,,
\label{e3_05}
\ee
Eq.~\eqref{e3_04} describes a {\em small} increment occurring from $n$th to $(n+1)th$ step:
\be
\F[\tv_{n+1}] = e^{i\beta(\omega^2-\opi^2) \dz} \, 
\F \left[ \tv_n + i\gamma\dz ( \usol^2\tv_n^* + 2|\usol|^2 \tv_n) \right]\,.
\label{e3_06}
\ee
Indeed, from \eqref{e3_06} and the estimate $\exp[i\beta(\omega^2-\opi^2) \dz]=1+O(\epsilon)$
it follows that $\tv_{n+1}-\tv_n=O(\epsilon)$.
Accordingly, we define a {\em continuous} variable $\tv(z,t)$ which at $z_n=n\dz$ is
related to $\tu_n$ by \eqref{e3_05}. Then, we show in Appendix 1 that a variable
\ $\tw=\tv\,e^{-iKz}$ \ satisfies an equation
%   %
%   \be
%   \big(\F[\tw]\big)_z=i\beta(\omega^2-\opi^2)\F[\tw]\,+\,\F\Big[-iK\tw + i\gamma U^2(\tw^*+2\tw)\Big]\,+
%    \, O(\epsilon),
%   \label{e3_06}
%   \ee
%   %
%
\be
\tw_z=-i\beta(\tw_{tt}+\opi^2\tw)-iK\tw+i\gamma U^2(\tw^*+2\tw),
\label{e3_07}
\ee
where $K$ and $U\equiv U(t)$ are defined in \eqref{e1_08}, and we have neglected terms of magnitude 
$O(\epsilon)$ (see Appendix 1). 
Note that the
first group of terms on the right-hand side of \eqref{e3_07} has the order of magnitude $O(1/\epsilon)$
(see the text after \eqref{e3_04}), and thus it
describes waves rapidly oscillating in both $z$ and $t$, as we have announced in Section 3.1.

We seek a solution of \eqref{e3_07} as
\be
\tw=p\,e^{-i\opi t} + m^*\,e^{i\opi t}\,,
\label{e3_08}
\ee
where $p(t,z)$ and $m(t,z)$ are assumed to vary in time on a scale $O(1)$. This agrees with
the numerical observation in Section 2 that frequencies of unstable Fourier modes differ
from $\pm\opi$ by $O(1)$ in most cases. 
Substituting \eqref{e3_08} into \eqref{e3_07}
and separating terms proportional to $\exp(\pm i\opi t)$, one obtains:
\bsube
\be
p_z=-2(\beta/\epsilon) p_t -iKp-i\beta p_{tt} + i\gamma U^2 (m + 2p), 
\label{e3_09a}
\ee
\be
m_z=2(\beta/\epsilon) m_t + iKm + i\beta m_{tt} - i\gamma U^2 (p + 2m),
\label{e3_09b}
\ee
\label{e3_09}
\esube
where we have used definition \eqref{e3_02}.
Using again the aforementioned numerical observation from Section 2, 
we seek solutions of \eqref{e3_09} in the form
\bsube
\be
p=p_{\rm slow}(\tau)\,\exp\left[ -i\Omega t + 2i(\beta/\epsilon)\Omega z + \beta\Lambda z \right],
\label{e3_10a}
\ee
\be
m=m_{\rm slow}(\tau)\,\exp\left[ i\Omega t + 2i(\beta/\epsilon)\Omega z + \beta\Lambda z \right],
\label{e3_10b}
\ee
\label{e3_10}
\esube
where 
\be
\Omega=O(1)
\label{add3_01}
\ee
labels a particular Fourier mode of $\tw$, parameter $\Lambda$ is to be determined later, 
and
\be
\tau=\epsilon t,  \qquad \Lambda \equiv \Lamr+i\Lami\,.
\label{e3_11}
\ee
Note that $\beta\Lamr$ equals Re$(\lambda)$ in \eqref{e1_05b}.
When writing that $p_{\rm slow}$ and $m_{\rm slow}$, as well as $P(\tau)$ and $M(\tau)$ in \eqref{e3_14}
below, are functions of the ``slow" time $\tau$, we mean that
\be
(p_{\rm slow})_t, \; (m_{\rm slow})_t, \; P_t, \; M_t \quad \mbox{are all of order $O(\epsilon)$}.
\label{e3_12}
\ee
The order-of-magnitude estimate \eqref{add3_01}, which we empirically deduced from numerical experiments,
will be generalized at the end of Section 4.

Substitution of \eqref{e3_10} into \eqref{e3_09} produces a pair of $z$-independent equations:
\bsube
\be
2(1+\epsilon\Omega)(p_{\rm slow})_{\tau} = \left( -iK/\beta +i\Omega^2 -\Lambda \right) p_{\rm slow}
  + i(\gamma/\beta)U^2(\tau/\epsilon)\, \left( m_{\rm slow}\,e^{2i\Omega \tau/\epsilon} + 2p_{\rm slow} \right),
\label{e3_13a}
\ee
\be
2(1-\epsilon\Omega)(m_{\rm slow})_{\tau} = \left( -iK/\beta +i\Omega^2 +\Lambda \right) m_{\rm slow}
  + i(\gamma/\beta)U^2(\tau/\epsilon)\, \left( p_{\rm slow}\,e^{-2i\Omega \tau/\epsilon} + 2m_{\rm slow} \right).
\label{e3_13b}
\ee
\label{e3_13}
\esube
In writing these equations we have neglected terms $O(\epsilon^2)$, which will be justified by 
subsequent calculations. Note also that the soliton background, $U^2(\tau/\epsilon)$, presents
a {\em narrow} obstacle for $p_{\rm slow}(\tau)$ and $m_{\rm slow}(\tau)$. 
Since outside the
soliton, $p_{\rm slow}$ and $m_{\rm slow}$ are not coupled to each other, 
we use yet another substitution:
\bsube
\be
p_{\rm slow} = P\,\exp \left[ \frac{ i(-K/\beta + \Omega^2) - \Lambda}{2(1+\epsilon\Omega)}\,\tau \right],
\label{e3_14a}
\ee
\be
m_{\rm slow} = M\,\exp \left[ \frac{ i(-K/\beta + \Omega^2) + \Lambda}{2(1-\epsilon\Omega)}\,\tau \right].
\label{e3_14b}
\ee
\label{e3_14}
\esube
Then $P$ and $M$ change only in the vicinity of the soliton according to
\bsube
\be
P_\tau=\frac{i\gamma\,U^2(\tau/\epsilon)}{2\beta (1+\epsilon\Omega)}\,
 \left( M\,e^{2i\Omega \cdot \tau/\epsilon + O(1)\cdot\tau} + 2P \right),
 \label{e3_15a}
\ee
\be
M_\tau=\frac{i\gamma\,U^2(\tau/\epsilon)}{2\beta (1-\epsilon\Omega)}\,
 \left( P\,e^{-2i\Omega \cdot \tau/\epsilon + O(1)\cdot\tau} + 2M \right).
 \label{e3_15b}
\ee
\label{e3_15}
\esube
Integrating these equations over the entire time window $-\epsilon T/2 \le \tau \le \epsilon T/2$,
one obtains:
\bsube
\be
P(+\epsilon T/2) - P(-\epsilon T/2) = i\epsilon\,\frac{\gamma}{2\beta} 
     \left( \F[U^2](-2\Omega)\,M(0) + 2\F[U^2](0)\,P(0) \right),
\label{e3_16a}
\ee
\be
M(+\epsilon T/2) - M(-\epsilon T/2) = i\epsilon\,\frac{\gamma}{2\beta} 
     \left( \F[U^2](2\Omega)\,P(0) + 2\F[U^2](0)\,M(0) \right).
\label{e3_16b}
\ee
\label{e3_16}
\esube
Here we have neglected terms $O(\epsilon^2)$ and used the definition \eqref{e1_03} of the
Fourier transform and the fact that the soliton is centered at $\tau=0$. 
These relations prompt us to denote relative jumps $\epsilon J_P$ and $\epsilon J_M$ that 
$P$ and $M$ undergo across the soliton:
\be
P(\pm \epsilon T/2) \equiv P(0)\,(1\pm \epsilon J_P/2), \qquad
M(\pm \epsilon T/2) \equiv M(0)\,(1\pm \epsilon J_M/2).
\label{e3_17}
\ee
Thus, by definition, \ $J_{P,M}=O(1)$. 
According to \eqref{e3_16} and with the same accuracy, these jumps satisfy
\bsube
\be
J_P = i\,\frac{\gamma}{2\beta} 
     \left( \F[U^2](-2\Omega)\,R + 2\F[U^2](0) \right),
\label{e3_18a}
\ee
\be
J_M = i\,\frac{\gamma}{2\beta} 
     \left( \F[U^2](2\Omega)\,\frac1R + 2\F[U^2](0) \right),
\label{e3_18b}
\ee
\label{e3_18}
\esube
where \ $R=M(0)/P(0)$.

Let us pause for a moment and recall what we are trying to do: We want to determine the instability
increment $\beta\Lamr$ for the Fourier mode whose frequency is related to $\Omega$ by \eqref{e3_10} and
\eqref{e3_08}. The two equations \eqref{e3_18} for three unknowns $J_{P,M}$ and $R$ are insufficient
for this purpose; note that they do not even involve $\Lambda$. The missing relations that will allow
us to complete our task are supplied by the periodicity condition satisfied by $P(\tau)$ and $M(\tau)$
and are obtained as follows. First, note that the numerical error
$\tu_n(t)$, and hence $\tw(t,z)$ defined before \eqref{e3_06}, satisfies periodic  boundary
conditions in $t$ by virtue of the split-step method \eqref{e1_02} using the discrete Fourier transform.
Second, since $p(t,z)$ and $m^*(t,z)$ have different $z$-dependences (see \eqref{e3_10}), then by virtue
of \eqref{e3_08} {\em each one} of these functions must satisfy periodic boundary conditions in $t$.
Third, along with \eqref{e3_08}, \eqref{e3_10}, and \eqref{e3_14}, this implies that
\bsube
\be
\frac{P(+\epsilon T/2)}{P(-\epsilon T/2)} = \exp \left[ 
 \frac{ i \left( K/\beta - \Omega^2 \right) +\Lambda }{2(1+\epsilon\Omega)} \, \epsilon T + 
 i(\opi+\Omega)T\, \right],
\label{e3_19a}
\ee
\be
\frac{M(+\epsilon T/2)}{M(-\epsilon T/2)} = \exp \left[ 
 \frac{ i \left( K/\beta - \Omega^2 \right) -\Lambda}{2(1-\epsilon\Omega)} \, \epsilon T 
  + i(\opi-\Omega)T\,  \right].
\label{e3_19b}
\ee
\label{e3_19}
\esube
Finally, using the identity $\exp[2i\pi N]=1$ for integer $N$ and Eqs. \eqref{e3_19} and \eqref{e3_17},
and neglecting terms $O(\epsilon^2)$, we obtain:
\bsube
\be
\left[ i \left( K/\beta - \Omega^2 \right) +\Lambda \right] (1-\epsilon\Omega) \, \epsilon T/2 + 
 i(\dopi+\Omega)T\, = \, 2i\pi N_P + \epsilon J_P,
 \label{e3_20a}
 \ee
\be
\left[ i \left( K/\beta - \Omega^2 \right) -\Lambda \right] (1+\epsilon\Omega) \, \epsilon T/2 + 
 i(\dopi-\Omega)T\, = \, 2i\pi N_M + \epsilon J_M,
 \label{e3_20b}
 \ee
\label{e3_20}
\esube
where $N_{P,M}$ are some integer numbers and $\dopi$ is ``the fractional part" of $\opi$: If 
$\opi=2\pi n/T$ (where $n$ is not necessarily an integer) and $N_{\pi}$ is the integer part of
$n$, then $\dopi \equiv 2\pi(n-N_{\pi})/T$. Note that \ $\dopi=O(1/T)$.

We can now justify neglecting the terms $O(\epsilon^2)$ in
\eqref{e3_13} and in subsequent calculations. Indeed, if such terms had been retained, they would
have contributed amounts $O(\epsilon^2)$ and $O(\epsilon^2\Omega\cdot\epsilon T)$ \ to \eqref{e3_20}.
The former amount would be a higher-order contribution
than that provided by terms $\epsilon J_{P,M}$, which we need to determine. The latter amount,
strictly speaking, depends on the order of magnitude of $(\epsilon T)$, but it will be clear
from our subsequent calculations that even terms \ $O(\epsilon\Omega\cdot\epsilon T)$ can be
neglected in the leading-order analysis.

We will now use Eqs. \eqref{e3_18} and \eqref{e3_20} to
determine for what values of $\Omega$ the instability increment $\beta\Lamr$ is nonzero.
To this end, we first subtract Eqs. \eqref{e3_20} from each other and take the real part, obtaining:
\be
\Lamr = \frac{{\rm Re}(J_P-J_M)}{T}\,.
\label{e3_21}
\ee
Next, adding Eqs. \eqref{e3_20} one obtains:
\be
\epsilon (J_P+J_M) = \left( iK/\beta - i\Omega^2  - \Lambda\, \epsilon\Omega \right) \epsilon T 
  + 2i\dopi T - 2 i \pi ( N_P+N_M)\,.
 \label{e3_22}
 \ee
 Using Eq.~\eqref{e3_21} one notices that the real part of the the right-hand side of 
 \eqref{e3_22} is of order $O(\epsilon^2)$ and hence should be neglected. Thus we conclude
 that in the main order, $(J_P+J_M)$ is purely imaginary.
 For future use, we also display the result of taking the 
imaginary part of the difference of the two equations \eqref{e3_20}:
 \be
 {\rm Im}(J_P-J_M) = \left( - \left( K/\beta - \Omega^2 \right) \epsilon\Omega + 
  \Lami  \right)\,\epsilon T \,+\, 2\Omega T - 2\pi (N_P-N_M)\,.
 \label{e3_23}
 \ee

Since $\Lamr$ is proportional to the real part of $(J_P-J_M)$, we solve for the latter quantity
using Eqs. \eqref{e3_18}. To that end, we first solve for $R$ by adding these equations and then
substitute the answer in their difference, obtaining:
\be
J_P-J_M= \pm \frac{i\gamma}{\beta} 
 \sqrt{ \big( 2\F [U^2](0) + i(\beta/\gamma)(J_P+J_M) \big)^2 - \big| \F [U^2](2\Omega) \big|^2 }\,.
 \label{e3_24}
 \ee
 Now recall that, as noted after \eqref{e3_22}, $(J_P+J_M)$ is purely imaginary. 
 Also, $\F[U^2](0)$ is real. Then 
 the real part of the right-hand side of \eqref{e3_24} is nonzero when
\be
 - \big| \F [U^2](2\Omega) \big| \;\le\; 2\F [U^2](0) + i(\beta/\gamma)(J_P+J_M)  \;\le\;
  \big| \F [U^2](2\Omega) \big|\,.
 \label{e3_25}
 \ee
 Under this condition, one also has
 \be
 {\rm Im}(J_P-J_M)=0.
 \label{e3_26}
 \ee
 %

 %  Let us outline the program of the remaining calculation.
 %
 Thus, the instability increment, $\beta\Lamr$, is found from Eqs. \eqref{e3_21} and \eqref{e3_24}, where
 $(J_P+J_M)$ is determined from \eqref{e3_22}. The last two equations, in their turn, involve
 three yet undetermined quantities:  
 $\Omega$, which labels the frequency of
 an unstable Fourier mode, parameter $\Lami$ introduced in  \eqref{e3_10}, and $(N_P+N_M)$. 
 We will now show that within the accuracy adopted in our calculations,
 the former two quantities enter all equations in a unique combination
 $(\Omega + \epsilon \Lami/2)$, and hence the number of yet undetermined quantities reduces to two.
 Indeed, 
 upon substitution of \eqref{e3_14} into \eqref{e3_10}, we observe that {\em up to terms $O(\epsilon^2)$},
 both $t$- and $z$-dependences of $p$ and $m$ involve $\Omega$ and $\Lami$ only in the
 aforementioned combination  $(\Omega+\epsilon \Lami/2)$. 
 Next, by inspection of formulae \eqref{e3_22} and \eqref{e3_23}, one can easily
 see that, within the same accuracy, they also involve $\Omega$ and $\Lami$ only in that combination.
 This means that  one can set $\Lami=0$. 
 Then from 
 \eqref{e3_23} and \eqref{e3_26} one finds:
 \be
 \Omega = \pi(N_P-N_M)/T\,,
 \label{e3_27}
 \ee
 where we have discarded terms $O(\epsilon^2)$.
 Then \eqref{e3_22} is rewritten as
 \be
  i (J_P+J_M) = - \left( \frac{K}{\beta} + \frac{2\dopi}{\epsilon} \right) T 
   + \left( \frac{\pi^2 (N_P-N_M)^2}{T} +  \frac{2 \pi ( N_P+N_M)}{\epsilon} \right)\,.
  \label{e3_28}
 \ee
 Finally, one substitutes the last two equations into \eqref{e3_25} and determines those
 values of $(N_P \pm N_M)$  where the instability can occur.
 The instability increment is computed from \eqref{e3_21} and \eqref{e3_24},
 and the frequency of the unstable mode, from \eqref{e3_27}.
 Examples of such a calculation, producing the theoretical results shown in
 Figs.~\ref{fig_3} and \ref{fig_4}, are given in Section 4.2.

 In Appendices 2 and 3 we present related technical results. Namely, in
 Appendix 2 we show that if one seeks unstable modes not specifically near $\pm\opi$,
 as we did at the beginning of this subsection,
 but instead near an arbitrary pair of frequencies $\pm\omega_0$, one discovers
 that the instability can arise only near $\pm\opi$, \ $\pm\omega_{2\pi}$, etc.
 In Appendix 3 we modify the analysis of this subsection to apply to a monochromatic-wave background
 \eqref{e1_04} and thereby recover results that can be deduced from those obtained by
 Weideman and Herbst \cite{WH}.

 % ---------------------------------------------------------------
 \section{Validation of the theory}
 \setcounter{equation}{0}

 In the first subsection below, we will
 give qualitative explanations of the instability features 
 described in Section 2: the locations and widths of the instability peaks, and
 high sensitivity of the instability to the step size, the time window length $T$,
 and the shape of the background solution. In the second subsection, we will work out
 an example showing how Eqs. \eqref{e3_21}, \eqref{e3_24}, \eqref{e3_25}, \eqref{e3_28},
 and \eqref{e3_27} were used to compute the increment and frequency of the unstable
 Fourier modes reported in Figs.~\ref{fig_3} and \ref{fig_4}.
 We will conclude by generalizing the order-of-magnitude estimate \eqref{add3_01}
 for the unstable modes' separation frequency, $2\Omega$.

 % ..............................................................................
 
 \subsection{Qualitative explanation of instability features reported in Section 2}

 The results of Section 3.2 allow us to explain why the instability peaks,
 shown in Fig. \ref{fig_2}, are not reflectionally symmetric about $\omega=0$ (see
 the paragraph after Eq.~\eqref{e2_02}). Indeed, from Eqs. \eqref{e3_08}, \eqref{e3_10},
 and \eqref{e3_14}, one sees that the frequencies of two coupled unstable modes are
 \be
 \big( \opi + \Omega - \epsilon (-K/\beta + \Omega^2)/2 \big) \qquad \mbox{and} \qquad
 -\big( \opi - \Omega - \epsilon (-K/\beta + \Omega^2)/2 \big) \,.
 \label{e3_29}
 \ee
 Thus, given the sign difference of the second terms inside the parentheses above,
 the mode at $\omega^{(+)}_{\rm right}$ is coupled to the mode at 
 $\omega^{(-)}_{\rm right}$ and {\em not} to that at $\omega^{(-)}_{\rm left}$,
 as it would have been in the case of reflectional symmetry.

 Similarly, two other features of the instability spectra reported in Section 2 can also be
 explained. First, note from
  \eqref{e3_27} and \eqref{e3_28} that if a value $\Omega>0$ is found to correspond to
 an instability peak, then so is $-\Omega<0$. 
 This observation, along with relations 
 \be
 \omega^{(\pm)}_{\rm right} = \pm \big( \big[\opi - \epsilon (-K/\beta + \Omega^2)/2 \big] \pm \Omega \big) 
 \qquad \mbox{and} \qquad
 \omega^{(\pm)}_{\rm left} = \pm \big( \big[ \opi - \epsilon (-K/\beta + \Omega^2)/2 \big]
  \mp \Omega  \big),
 \label{e3_30}
 \ee
 which follow from \eqref{e3_29}, explains why relations \eqref{e2_03} hold.
 Second, the slight asymmetry of the frequencies $\omega^{(\pm)}_{\rm right}$ and
 $\omega^{(\pm)}_{\rm left}$ about the respective $\pm\opi$, is also easily explained.
 Indeed, the frequencies of, say, the peaks
 at $\omega^{(+)}_{\rm right}$ and $\omega^{(+)}_{\rm left}$, as seen from \eqref{e3_30}, are
 centered about \ $\big[\opi - \epsilon (-K/\beta + \Omega^2)/2 \big]$ and {\em not} about $\opi$.
 The maginitude of the shift, $ - \epsilon (-K/\beta + \Omega^2)/2$, agrees
 with the numerically observed values. For example, in the experiments reported in Table 1, 
 $K=1$, $\beta=-1$, $\epsilon \sim 0.04$, and $\Omega \sim 0.5$, and hence 
 $ -\epsilon (-K/\beta + \Omega^2)/2 \sim -0.025$. This should be compared to the
 experimental values of $\big[(\omega^{(+)}_{\rm left}-\opi)+(\omega^{(+)}_{\rm right}-\opi)\big]/2$,
 which from Table 1 are seen to vary between $-0.02$ and $-0.03$ \footnote{A more detailed comparison
 would require keeping at least one more significant digit in the data of Table 1, but
 such a comparison does not appear to be needed. Rather, it is the quantitative
 agreement of our theory and numerical experiments, reported in Figs.~\ref{fig_3} and \ref{fig_4}
 and presented in the next subsection, which seems to be the most important test confirming the
 validity of the theory.}.

 Next, we can explain why in most cases, as mentioned in Sections 2 and 3.1,  the
 instability peaks contain just one node\footnote{As we already mentioned in Section 2,
 the pedestals around the peaks arise due to a non-numerical --- modulational --- instability,
 and those pedestals are hence unrelated to the foregoing explanation.}. 
 Consider Eqs. \eqref{e3_27} and \eqref{e3_28}
 and assume that at frequency labeled by
 $\Omega_0$, corresponding to a particular $(N_P-N_M)_0$, there is an instability. 
 The frequency at the adjacent node differs from this $\Omega_0$ by $2\pi/T$, and hence,
 according to \eqref{e3_27}, 
 the corresponding $(N_P-N_M)$ differs from $(N_P-N_M)_0$ by 2. Consistently with this,
 one can take $(N_P+N_M)=(N_P+N_M)_0$. 
 Then the corresponding value of the right-hand side of \eqref{e3_28} differs 
 from that value at $\Omega_0$ by $2\pi^2(N_P-N_M)_0/T\,\cdot 2 = 4\pi\Omega_0$.
 On the other hand, the interval of values where $i(J_P+J_M)$ corresponds to an instability
 is found from \eqref{e3_25} to have the width \ $2\,\big|(\gamma/\beta) \F [U^2](2\Omega_0) \big|$.
 Thus, whenever
 \be
  4\pi\Omega \,>\, 2\,\big|(\gamma/\beta) \F [U^2](2\Omega) \big|,
  \label{e3_31}
  \ee
 the instability peak can contain only one node. Using expression \eqref{e1_08} for 
 $U(t)$, one finds that condition \eqref{e3_31} holds for $\Omega > 0.46$. In other words,
 it is only when the peaks at $\omega_{\rm left}$ and $\omega_{\rm right}$ are separated by less
 than $0.92$ that they can contain more than one node. In practice, the separation should be even smaller,
 given a random location of the node with the unstable frequency within the instability interval.
 In fact, among the simulations reported in Section 2, we observed multiple nodes per peak
 only when $(\omega_{\rm right}-\omega_{\rm left})$ was about 0.25 or less.

 The high sensitivity of the instability to the length $T$ of the time window, which was
 highlighted in Section 2 (see also Table 2 there), is also easily
 explained using \eqref{e3_28}. Suppose this length is changed from $T_0$ to $T$, so that the
 relative change \ $(T-T_0)/T_0$ is small. Then the
 right-hand side of \eqref{e3_28} is changed by an amount
 \be
 - \left( \frac{2\dopi}{\epsilon} + \Omega^2 \right) \cdot T_0 \cdot \frac{T-T_0}{T_0}\,.
 \label{e3_32}
 \ee
 The coefficient in front of $(T-T_0)/T_0$ is large due to $T_0$ being large and the terms in
 parentheses being $O(1)$. Note that this does not necessarily imply that the instability
 will disappear, since new values of $N_P$ and $N_M$ may exist for which the right-hand side of
 \eqref{e3_28} will be inside the instability interval \eqref{e3_25}; see the next subsection
 for a quantitative example in a similar situation where not $T$ but $\dz$, is varied.

 Finally, while instability on the background of several well-separated solitons cannot be 
 computed from Eqs.~\eqref{e3_24} and \eqref{e3_25}, which were obtained for a single
 soliton, one can still explain why the instability is sensitive to the number of the solitons and
 their relative locations. Indeed, the number of solitons determines the number of jumps in the
 ``slow" variables $P$ and $M$; see \eqref{e3_15}--\eqref{e3_18}. Its is those jumps that allow
 modes that exponentially grow/decay (in $t$) away from the jumps to exist in the presence of
 periodic boundary conditions in $t$; see \eqref{e3_20}. And its is those modes, exponentially growing
 or decaying (as opposed to purely oscillating) {\em in $t$}, that are also exponentially growing
 (i.e., are unstable), in the evolution variable $z$; note the same exponent $\Lambda$ in \eqref{e3_14}
 and in \eqref{e3_10}. As for the relative locations of the solitons, they (as well as the
 solitons' phases) determine the counterparts of parameter $R$ in a generalization of Eqs.~\eqref{e3_18}
 for a multi-soliton background.

 % ....................................................................
 
 \subsection{Example of calculation of instability increment and frequency}
 
 Here we will explain in detail how the instability increment and frequency for the generic case
 are found. Such a generic case corresponds to, e.g., $\dz=0.0043$ in Fig.~\ref{fig_3},
 for which we will present the calculations below. 
 Its instability spectrum is qualitatively the same as that shown in Fig.~\ref{fig_2}.
 Then we will comment on less common cases, whose spectra are shown in Fig.~\ref{fig_5}. 
 We will conclude this section by establishing a general 
 dependence of $\Omega$ as a function of $\epsilon$ and $T$.

 For all cases considered here, the parameters in the background
 soliton \eqref{e1_08} are: \ $\beta=-1$, \ $\gamma=2$, $A=1$, whence $K=1$ and
 $U(t)={\rm sech}(t)$. The time window
 is $T=128\pi$ and the number of grid points is $2^{12}$. 
 These are the same parameters as were used in obtaining Figs.~\ref{fig_3} -- \ref{fig_5}.

 To begin, let us recall from \eqref{e3_24} and \eqref{e3_25}
 that for the mode labeled by $\Omega$ to be unstable, the left-hand side
 of \eqref{e3_28} must fall within the interval
 \be
 \left( 4\F [U^2](0) -  2\big| \F [U^2](2\Omega) \big|\,,\;\;
  4\F [U^2](0) +  2\big| \F [U^2](2\Omega) \big|\, \right)\,
 \label{e4_05}
 \ee
 (here we have used that $\gamma/\beta=-2$).
 While $\Omega$ is yet undetermined and hence $\big| \F [U^2](2\Omega) \big|$ is not known
 exactly, it cannot exceed $\F [U^2](0)=2$, and hence interval \eqref{e4_05} is always inside
 the interval $(4,\;12)$. Thus, this interval contains values of order 1, i.e. much smaller
 than any of the terms on the right-hand side of \eqref{e3_28}, which are proportional to
 the large parameters $T$ and $1/\epsilon=\opi$. This observation motivates our strategy of
 finding suitable values of $(N_P\pm N_M)$.

 First,  given the above values of parameters, one finds $\opi=1/\epsilon\approx27.03$,
 $\dopi\approx0.0140$, and then the first term
 on the right-hand side of \eqref{e3_28}:
 \be
 - \left( \frac{K}{\beta} + \frac{2\dopi}{\epsilon} \right) T \approx 96.86\,.
 \label{e4_06}
 \ee
 Next, since the first term in the second parentheses, \ $(\pi(N_P-N_M))^2/T$, \ is always
 positive, we seek such an integer value for $(N_P+N_M)$ that an expression
 \be
 - \left( \frac{K}{\beta} + \frac{2\dopi}{\epsilon} \right) T  + 
 \frac{2 \pi ( N_P+N_M)}{\epsilon} 
 \label{e4_07}
 \ee
 is as close to zero as possible but negative\footnote{Modifications of this requirement
 will give rise to ``less common" cases, considered later.}. 
 Since $2\pi/\epsilon\approx 169.84$, the corresponding $N_P+N_M=-1$, yielding the value
 of $-72.98$ for \eqref{e4_07}. 
 Finally, an integer value of $(N_P-N_M)$ is sought to make the value of the
 right-hand side of \eqref{e3_25} fit within the interval \eqref{e4_05}. To this end, one
 first {\em estimates} $(N_P-N_M)$ by forcing that right-hand side to equal zero:
 \be
 (N_P-N_M)_{\rm estimated} = \sqrt{ \frac{T}{\pi^2} \left[
  \left( \frac{K}{\beta} + \frac{2\dopi}{\epsilon} \right) T  -
 \frac{2 \pi ( N_P+N_M)}{\epsilon}  \right] }\,.
 \label{e4_08}
 \ee
 Here $(N_P-N_M)_{\rm estimated}$ is not an integer; e.g., for the example above it is
 approximately $54.53$. Then a nearby integer value for $(N_P-N_M)$ is found by trial and error
 so that it makes the right-hand side of \eqref{e3_28} fit within the interval \eqref{e4_05}.
 For the above example, this value is $(N_P-N_M)=57$, yielding $\Omega\approx 0.44$ (see
 \eqref{e3_27})  and $\beta\Lamr=0.0067$ (see \eqref{e3_24} and \eqref{e3_21}). Comparison of
 so computed $\beta\Lamr$ and $\Omega$ with respective values found numerically is shown
 in Figs.~\ref{fig_3} and \ref{fig_4} for various values of $\dz$.

 Two notes are in order about the last step described above. First, since both $N_P$ and
 $N_M$ are integers, then $(N_P+N_M)$ and $(N_P-N_M)$ are either both odd or both even.
 Second, once a guess is made about $(N_P-N_M)$, then $\Omega$ is computed from \eqref{e3_27}
 and hence $\big| \F [U^2](2\Omega) \big|$ in \eqref{e4_05} and \eqref{e3_24} is 
 found from \eqref{e1_03}.

 We will now briefly describe what we have referred to above as ``less common" cases,
 depicted in Fig.~\ref{fig_5}. Let us first address the situation depicted in Fig.~\ref{fig_5}a,
 which occurs for, e.g., $\dz=0.0044,\,0.0045, \, 0.00474,\, 0.0049$. The inner pair of instability peaks
 is found in these cases as described above. For example, for $\dz=0.0049$, it is obtained using
 $(N_P+N_M)=-2$ and $(N_P-N_M)=62$. The corresponding calculated instability increment 
 and the peak separation are $\beta\Lamr\approx 0.0051$ and $2\Omega\approx2\cdot 0.48$, which are very close
 to the experimentally observed values. 
 However, if for the same value of $\dz$ one takes $(N_P+N_M)=-3$,
 then one finds that an instability can exist at $\Omega\approx 0.78$, corresponding to $(N_P-N_M)=101$.
 The calculated and observed values of the instability increment for this outer pair of peaks are
 $0.0030$ and $0.0039$, respectively.

 Let us note that we checked --- by trial and error, as described in the previous paragraph ---
 for the existence of a secondary pair of instability peaks for every value of $\dz$ reported in
 Figs.~\ref{fig_3} and \ref{fig_4}. For each $\dz$ where such a secondary pair was observed in numerics,
 we also found it analytically, with the agreement between the calculated and observed values of $\Omega$
 being excellent and those of $\Lamr$ being good (the discrepancy between such values for $\dz=0.0049$,
 reported above, was the worst one we found). 
 We even found, both numerically and analytically, small tertiary peaks for $\dz=0.0044$. 
 On the other hand, we {\em analytically} found secondary peaks
 for $\dz=0.0042,\, 0.00475, \,0.00477$ where they were not observed in numerics at $z=z_{\rm max}=2000$.
 However, at smaller $z$, such secondary peaks were indeed observed. The reason they were not
 observed for $z=2000$ was that they were ``drowned" by the pedestal of unstable modes which ``rose" around the
 primary peaks due to modulational --- i.e., non-numerical --- instability about the primary peaks.
 
 We will now comment on the calculation of the central peak for the case of $\dz=0.0050$,
 shown in Fig.~\ref{fig_5}b. Here the value of \eqref{e4_07} with $(N_P+N_M)=-1$ is approximately $11.1$.
 This is not negative, as we wanted this expression to be in the previous calculations (see the line
 after \eqref{e4_07}). However, for $(N_P-N_M)=0$, the right-hand side of \eqref{e3_28} is 
 $11.1\in(4,\,12)$, i.e. it is within the instability interval \eqref{e4_05} with $\Omega=0$ (see \eqref{e3_27}).
 According to the discussion found around Eq. \eqref{e3_31}, the instability peak at $\Omega=0$ should contain
 several grid nodes, which is confirmed by Fig.~\ref{fig_5}b.

 Finally, let us revisit the estimate for the order of magnitude of $\Omega$, which determines
 the frequency separation of the primary instability peaks. In \eqref{add3_01} we stated, as an empirical
 observation, that $\Omega=O(1)$. We will now use Eq.~\eqref{e4_08} to generalize it.
 Note that for the largest integer value of $(N_P+N_M)$ that still renders expression \eqref{e4_07}
 negative, the expression in the square brackets in \eqref{e4_08} is no greater than \ $2\pi\cdot 1/\epsilon$.
 The corresponding $(N_P-N_M)_{\rm estimated}$, and hence $(N_P-N_M)$, is then of the order 
 $O\big(\sqrt{T/\epsilon}\big)$. Substitution of this into \eqref{e3_27} yields
 \be
 \Omega = O \big( 1/\sqrt{\epsilon T} \big).
 \label{e4_09}
 \ee
 Thus, as one makes the time window wider in the order-of-magnitude sense 
 while keeping the other parameters fixed, the separation
 between the primary unstable Fourier modes should, {\em on average}, descrease. An implication of this
 for observing the numerical instability is presented in Section 5.

 % --------------------------------------------------
 \section{Conclusions}
 \setcounter{equation}{0}
 
 We reported, and then analytically explained, a numerical instability in the split-step Fourier method
 \eqref{e1_02} applied to the nonlinear Schr\"odinger equation \eqref{e1_01} with
 the background solution being the soliton \eqref{e1_08}. Properties of this instability, such as
 the dependence of its increment and unstable mode's location on the step size, numerical domain's
 length $T$, and details of the background solution, are quite different from those properties
 on the background of the monochromatic wave  \eqref{e1_04}, previously
 obtained by the von Neumann analysis 
 in \cite{WH}. Namely, the dependence of the instability on those parameters is seemingly very
 irregular, as we illustrated in Tables 1, 2 and Figs.~\ref{fig_3}, \ref{fig_4}. This is to
 be contrasted with monotonic dependences of the numerical instability observed for 
 %   the same  scheme \eqref{e1_02} on a monochromatic-wave background \cite{WH} or for 
 finite-difference
 schemes in both constant- and variable-coefficients equations, as described in textbooks
 on numerical methods. Moreover, we demonstrated (see the end of Section 2)
 that the principle of ``frozen coefficients" is not valid for the split-step Fourier method on
 the background of a localized nonlinear wave. In particular, the instabilities on the background
 of a single soliton, on one hand, and several well-separated identical solitons,
 on the other, can be drastically different.

 The analysis presented in Section 3 (see also Appendix 2)
 revealed that unstable modes can be found only near the resonant frequencies 
 $\pm\opi$, $\pm\omega_{2\pi}$, etc. (see \eqref{e1_06} and the sentence after that equation). 
 Interestingly, far away from the soliton, these unstable modes are {\em not} exactly periodic in the
 spatial variable $t$. Rather, their spatial envelope is exponentially growing or decaying
 in $t$ (see Eqs.~\eqref{e3_08}, \eqref{e3_10}, \eqref{e3_14}, and \eqref{e3_21}). 
 What makes these modes satisfy the periodic boundary conditions, which are implicitly
 imposed by the use of Fourier transform in \eqref{e1_02}, is the change that they undergo
 near the soliton. Thus, the finite size of the numerical domain $T$ is critical in our
 instability analysis, which is in stark contrast with the standard von Neumann analysis.

 As we noted above, the dependence of the instability increment on the parameters of the
 soliton and the numerical scheme is irregular, and there is no means to predict it ``quickly",
 i.e., bypassing the procedure illustrated in Section 4.2. 
 Again, this is in contrast with the instability analysis on a
 monochromatic-wave background \cite{WH} (see \cite{J} and Appendix 3). 
 However, one can still obtain two general conclusions: \ (i) what maximum value this increment 
 can have and \ (ii) how the chances to observe a numerical instability depend on the size $T$ of the
 numerical domain. 
 
 The maximum instability increment is obtained from Eqs. \eqref{e3_21} and \eqref{e3_24} by setting 
  \ $2\F [U^2](0) + i(\beta/\gamma)(J_P+J_M)=0$. This yields:
  \be
  \big(\beta\Lamr\big)_{\rm max} \,=\, 
  \frac1T \big( \gamma \big| {\mathcal F} [U^2](2\Omega) \big| \big)_{\rm max} \,=\, 
  \frac{\gamma}T {\mathcal F} [U^2](0) \equiv \frac{\gamma}T \int_{-T/2}^{T/2} U^2(t)\,dt\,.
  \label{e5_01}
  \ee
  Interestingly, the last expression above is the same as an analogous expression in the 
  monochromatic-background case (see Eq. \eqref{a3_05} in Appendix 3), where $U(t)\equiv A$.

 Now let us show that as the time window length $T$ is {\em substantially} decreased, the
 chances to observe numerical instability in any given simulation using \eqref{e1_02} with 
 \eqref{e2_01}, also decrease\footnote{Note the word `substantially'. As illustrated by Table 2
 in Section 2, changing $T$ by {\em only a fraction of its value} affects the occurrence
 of the instability in a non-monotonic and seemingly irregular way.}.
 This conclusion follows from a combination of arguments that led to formulae \eqref{e3_31}
 and \eqref{e4_09}. Indeed, recall from a discussion after Eq.~\eqref{e4_08} that an instability
 would arise only if near the non-integer number $(N_P-N_M)_{\rm estimated}$, there is an integer
 number $(N_P-N_M)$ that would make the right-hand side of \eqref{e3_28} fit within the interval
 \eqref{e4_05}. A sufficient condition that would guarantee that such an $(N_P-N_M)$ can be
 found is obtained similarly to \eqref{e3_31}: 
  \be
   4\pi\Omega_{\rm estimated} \,<\, 2\,\big|(\gamma/\beta) \F [U^2](2\Omega_{\rm estimated}) \big|,
   \label{e5_02}
   \ee
  where \ $\Omega_{\rm estimated} \equiv \pi (N_P-N_M)_{\rm estimated}/T$.
 As follows from the discussion around \eqref{e4_09}, \ 
 $\Omega_{\rm estimated}=O(1/\sqrt{\epsilon T})$. Thus, as $T$ decreases, the chances that
 condition \eqref{e5_02} may be satisfied, also decrease. Hence the smaller $T$ is, the ``less often"
 a numerical instability of scheme \eqref{e1_02} on the background of a soliton would be observed.
 On the other hand, the smaller $T$ is, the stronger the numerical instability, if it {\em is} observed,
 is on average; this follows from \eqref{e5_01}.
 Both these conclusions agree with our observations in Section 2: Compare Table 1, obtained for $T=32\pi$, with
 Fig.~\ref{fig_3}, obtained for $T=128\pi$.

 %  -------------------------------------------------------------------------------
 \section*{Acknowledgement}
 
 I thank Jianke Yang for drawing my attention to the problem considered in this paper and for
 stimulating discussions at the early stage of this work.

 %  -------------------------------------------------------------------------------

  \setcounter{equation}{0}
  \renewcommand{\theequation}{A1.\arabic{equation}}
\section*{Appendix 1: \ Derivation of Eq.~\eqref{e3_07}}
 
 We verified that one can derive \eqref{e3_07} from \eqref{e3_06} by using a Taylor series expansion
 of $\exp[i\beta(\omega^2-\opi^2)\dz]$ and $(\tv_{n+1}-\tv_n)$ in powers of $\epsilon$.
 However,  an alternative derivation presented below is much less tedious and,
 importantly, more intuitive.

 Let us first note that the last equation in \eqref{e1_02} is equivalent to the nonlinear Schr\"odinger
 equation \eqref{e1_01} {\em plus} a term proportional to 
 \be
 \dz\beta\gamma\cdot[\partial_{tt},|u(t,z)|^2]u(t,z) + O(\dz^2),
 \label{a1_01}
 \ee
 where $[ \ldots\,,\ldots]$ denotes a commutator. This follows from the Baker--Campbell--Hausdorff
 formula; see, e.g., Sec. 2.4.1 in \cite{Agrawal_book}. Next, Eq.~\eqref{e3_03} is a linearized
 version of the last equation in \eqref{e1_02}. Therefore, it must be equivalent to the linearized 
 nonlinear Schr\"odinger equation plus terms
 of order $O(\beta\dz\partial_{tt})=O(\beta\omega^2\dz)=O(\epsilon^2)$, 
 provided that we assume that $\omega \sim \partial_t = O(1)$,
 or, equivalently, that the central frequency in expanding $\exp[i\beta\omega^2\dz]$ in a 
 Taylor series is $0$.

 In writing \eqref{e3_04} and then \eqref{e3_06},
 we stated that the central frequency is $\opi$ (or $-\opi$) rather than $0$. 
 Correspondingly, Eq.~\eqref{e3_06} must be equivalent to a {\em modified} linearized 
 nonlinear Schr\"odinger equation written for a
 small deviation $\tv$, plus terms of order $O\big(\beta(\omega^2-\opi^2)\dz\big)=O(\epsilon)$
 (see the text after \eqref{e3_04}). 
 Here the modification consists in replacing the operator $\partial_{tt}$, 
 whose Fourier symbol is \ 
 $-\omega^2\equiv -(\omega^2-0^2)$, \  with the operator $\partial_{tt}+\opi^2$, 
 whose Fourier symbol is  \ $ -(\omega^2-\opi^2)$. Thus, \eqref{e3_06} in the time domain is
 \be
 \tv_z = -i\beta(\tv_{tt}+\opi^2\tv) + i\gamma (\usol^2 \tv^* + 2|\usol|^2 \tv) +O(\epsilon)\,.
 \label{a1_02}
 \ee
 Substituting into this equation $\usol$ from \eqref{e1_08}, changing the variable
 $\tv=\tw\exp[iKz]$, and neglecting the $O(\epsilon)$ term, one obtains Eq.~\eqref{e3_07}.

  %  -------------------------------------------------------------------------------
   
 \setcounter{equation}{0}
 \renewcommand{\theequation}{A2.\arabic{equation}}
\section*{Appendix 2: \ Location of instability peaks}

Here we will present an explaination
of why the frequencies of unstable modes must be near $\pm\opi$, $\pm \omega_{2\pi}$, etc. 
 
 Suppose we seek the instability near a pair of frequencies $\pm\omega_0$; i.e., we assume that
 \be
 |\omega - \omega_0| = O(1) \qquad \mbox{or} \qquad |\omega - (-\omega_0)| = O(1)\,.
 \label{a2_00}
 \ee
 Then, proceeding as explained in the
 text after Eq.~\eqref{e3_03}, we obtain an equation similar to \eqref{e3_04}, where
 the ``$(-1)$" and $\opi$ on the right-hand side are replaced with \ $\exp[i\phi_0]$ and $\omega_0$,
 respectively, where
 \be
 \phi_0=\beta\omega_0^2\dz\,.
 \label{a2_01}
 \ee
 Then \eqref{e3_05} and \eqref{e3_06} get replaced with 
 \be
 \tv_n=e^{-i\phi_0 n}\,\tu_n\,,
 \label{a2_02}
 \ee
 \be
 \F[\tv_{n+1}] = e^{i\beta(\omega^2-\omega_0^2) \dz} \, 
 \F \left[ \tv_n + i\gamma\dz ( \usol^2\tv_n^*\cdot\underline{e^{-2i\phi_0 n}} + 2|\usol|^2 \tv_n) \right]\,,
 \label{a2_03}
 \ee
where the term making the key difference between \eqref{e3_06} and \eqref{a2_03} is underlined.
Let us now note that if the phase rotation, \ $-2\phi_0$, \ 
in that term would equal $-2\pi N$, where $N$ is any integer,
then that term would equal $1$, and the subsequent analysis would proceed as in Section 3.2 without any changes.
Therefore, we can say that the {\em nontrivial} phase rotation in \eqref{a2_03} is 
\ $-(2\phi_0-2\pi N_0)$, where $N_0$ is the nearest integer to $\phi_0/\pi$. For simplicity, but
without loss of generality, let us assume that $N_0=1$; the case of $N_0\neq 1$ is completely
analogous. Then, Eq.~\eqref{a2_03} becomes
\be
 \F[\tv_{n+1}] = e^{i\beta(\omega^2-\omega_0^2) \dz} \, 
 \F \left[ \tv_n + i\gamma\dz \big( \usol^2\tv_n^*\, e^{-2i\beta(\omega_0^2-\opi^2)n\dz } + 
 2|\usol|^2 \tv_n \big) \right]\,,
 \label{a2_04}
\ee
where in rewriting the exponential term we have used \eqref{a2_01} and \eqref{e1_06}.

To go from the  discrete equation \eqref{a2_04}
to a counterpart of the continuous equation \eqref{e3_07},
we make two observations.
First, $n\dz=z$ in the second exponential term in \eqref{a2_04}.
Second, and perhaps counter-intuitively: \ 
Despite the presence of this possibly fast-oscillating exponential, Eq.~\eqref{a2_04} 
still describes a {\em small} change for $(\tv_{n+1}-\tv_n)$. This is due to the presence of
the small terms $i\gamma\dz$ and $\beta(\omega^2-\omega_0^2)\dz$ (see \eqref{a2_00}) on the
right-hand side of that equation. 
Therefore, the continuous variable $\tv(t,z)$ interpolating the discrete variable in \eqref{a2_04}
satisfies a counterpart of \eqref{a1_02}:
 \be
 \tv_z = -i\beta(\tv_{tt}+\omega_0^2\tv) + i\gamma \big(\usol^2 \tv^* \, e^{-2i\beta(\omega_0^2-\opi^2)z }
 + 2|\usol|^2 \tv \big) +O(\epsilon)\,.
 \label{a2_05}
 \ee
 In analogy with the argument presented in
Appendix 1, we make a change of variables
\be
\tv=\tw\,\exp\big[ i(K-\beta(\omega_0^2-\opi^2))z \big]\,,
\label{a2_06}
\ee
which transforms \eqref{a2_05} into the following counterpart of \eqref{e3_07}:
\be
\tw_z=-i\beta(\tw_{tt}+\omega_0^2\tw)-i\big(K - \beta(\omega_0^2-\opi^2)\big) \tw+i\gamma U^2(\tw^*+2\tw)\,.
\label{a2_07}
\ee
Then, a substitution analogous to \eqref{e3_08} with $\opi$ being replaced by $\omega_0$ 
into Eq.~\eqref{a2_07} yields a system of equations that is similar to \eqref{e3_09},
with the only changes being the replacements:
\be
\opi \;\; \mbox{by} \;\; \omega_0 \qquad \mbox{and} \qquad K \;\; \mbox{by}
\;\; (K- \beta(\omega_0^2-\opi^2))\,.
\label{a2_08}
\ee

We now need to consider two cases: \ (i) \ $\beta(\omega_0^2-\opi^2)=O(1)$ \ and \ (ii) \
$|\beta(\omega_0^2-\opi^2)| \gg 1$. We will show that in the first case, the results of
analysis of Eqs.~\eqref{e3_09} with replacements
\eqref{a2_08} reduce to those obtained in Section 3.2, and in the
second case, no instability can arise.

In case (i), \ $\omega_0-\opi=O(\epsilon)$, \ 
i.e. this case differs from that considered in
Section 3.2 only by a slight shift of the central frequency. 
Intuitively, such a shift could not change the location of the unstable peaks which we
found to be away from $\opi$ by an amount of approximately $\Omega=O(1)$. 
Formally, this can be justified by a tedious calculation that reveals that
Eqs.~\eqref{e3_27}, \eqref{e3_28} with 
replacements \eqref{a2_08} 
yield the same $\Omega$ as the original Eqs.~\eqref{e3_27}, \eqref{e3_28}.
Then, Eqs.~\eqref{e3_10}, \eqref{e3_14} with 
replacements \eqref{a2_08}
yield the same $t$-dependence of the solution of \eqref{e3_08} as the
original Eqs.~\eqref{e3_10}, \eqref{e3_14}.
Thus, in case (i), the parameters of the instability
reduce to those found in Section 3.2.

In case (ii), one cannot proceed as in case (i) by merely using replacements \eqref{a2_08}
in Eqs.~\eqref{e3_09}. The reason is that $\big| K-\beta(\omega_0^2-\opi^2)\big| \gg 1$,
whereas in case (i) one had $\big( K-\beta(\omega_0^2-\opi^2)\big) =O(1)$.
Indeed, in case (ii), a substitution \eqref{e3_10} with $\Omega=O(1)$ (which is our starting assumption ---
see \eqref{a2_00}) would not yield $p_{\rm slow}$ and $m_{\rm slow}$ that
would be slow functions of $t$; see the first term on the right-hand side of \eqref{e3_13}. 
The only way the large term $\big( K-\beta(\omega_0^2-\opi^2)\big)$ could be eliminated
from the couterpart of \eqref{e3_09} is by using {\em different} $z$-dependences in the
exponentials in \eqref{e3_10}:
\bsube
\be
p=p_{\rm slow}(\tau,z)\,\exp\left[ -i\Omega t + 2i(\beta/\epsilon)\Omega z + \beta\Lambda z 
 -i  \big( K-\beta(\omega_0^2-\opi^2)\big)z  \right],
\label{a2_09a}
\ee
\be
m=m_{\rm slow}(\tau,z)\,\exp\left[ i\Omega t + 2i(\beta/\epsilon)\Omega z + \beta\Lambda z 
 + i  \big( K-\beta(\omega_0^2-\opi^2)\big)z \right]\,.
\label{a2_09b}
\ee
\label{a2_09}
\esube
(Note that the last terms in the exponents in \eqref{a2_09} essentially undo 
transformation \eqref{a2_06}.) 
In \eqref{a2_09}, $p_{\rm slow}$ and $m_{\rm slow}$ are slow functions of $t$, but not of $z$.
Substituting \eqref{a2_09} into the counterpart of
\eqref{e3_09} one would obtain, instead of the $z$-independent system \eqref{e3_13},
a $z$-dependent system of the form:
\bsube
\be
(p_{\rm slow})_z = b_1 (p_{\rm slow})_{\tau} + b_2 p_{\rm slow}
  + b_3\,  m_{\rm slow}\,e^{ 2i \big( K-\beta(\omega_0^2-\opi^2)\big)z },
\label{a2_10a}
\ee
\be
(m_{\rm slow})_z = c_1 (m_{\rm slow})_{\tau} + c_2 m_{\rm slow}
  + c_3\, p_{\rm slow}\,e^{ - 2i \big( K-\beta(\omega_0^2-\opi^2)\big)z },
\label{a2_10b}
\ee
\label{a2_10}
\esube
where all the coefficients $b_1$ through $c_3$ are of order $O(1)$ and independent of $z$.
The presence of rapidly oscillating exponential terms in \eqref{a2_10} makes the effect
of the coupling terms negligible, and then system \eqref{a2_10} gets essentially decoupled
into two independent equations for $p_{\rm slow}$ and $m_{\rm slow}$, which does not
exhibit any instability. Thus, in case (ii) numerical instability does not occur.

 %  -------------------------------------------------------------------------------
   
\setcounter{equation}{0}
\renewcommand{\theequation}{A3.\arabic{equation}}
\section*{Appendix 3: \ Instability on the background of a monochromatic wave}

Here we will use the method presented in Section 3.2 to find the location and growth rate of
the numerically unstable Fourier modes of method \eqref{e1_02} on the background of a monochromatic wave
\eqref{e1_04} with $\Omega_{\rm cw}=0$.  Let us note that these results can be obtained from
formulae (65), (37), and (64) of \cite{WH} by expanding them in a power series of the step size $\dz$
(denoted there by $\tau$).
In such a way, the growth rate of the most unstable mode was obtained in \cite{J}.

The starting point of our derivation is system \eqref{e3_09}, which holds both for the soliton and
monochromatic-wave backgrounds. In the latter case, $U(t)\equiv A$ and $K=\gamma A^2$, where without
loss of generality we assume that $A$ is real. Thus, now
this system, unlike \eqref{e3_09} on the soliton background, has {\em all constant} coefficients,
and hence we can look for its solution in the form
\be
\{p,\,m\} \,=\, \{P,\,M\}\,e^{iW \epsilon t +  \lambda z}\,, \qquad W=O(1).
\label{a3_01}
\ee
Note that unlike in \eqref{e3_10}, here the $p$- and $m$-components of the small deviation $\tilde{w}$
have the same $t$-dependence. Also, we have used the notation $W$, not $\Omega$, in \eqref{a3_01}, because,
unlike $\Omega$, the variable $W$ does not have the dimension of frequency.  Rather,
$\epsilon W \equiv W/\opi$ has the same dimension as $\Omega$.

Substitution of \eqref{a3_01} into \eqref{e3_09} with the aforementioned values of $U$ and $K$ yields:
\bsube
\be
(2i\beta W-i\gamma A^2 +\lambda)\,P - i\gamma A^2\, M = 0,
\label{a3_02a}
\ee
\be
-i\gamma A^2 \,P + (2i\beta W-i\gamma A^2 -\lambda)\, M = 0,
\label{a3_02b}
\ee
\label{a3_02}
\esube
where we have neglected terms $O(\epsilon^2)$. Then the instability growth rate is
\be
\lambda= \sqrt{ (\gamma A^2)^2 - (2\beta W - \gamma A^2)^2 }\,.
\label{a3_03}
\ee
The location of the unstable mode(s) follows from the definition of the mode's frequency, 
$\omega=\opi-\epsilon W \equiv \opi - (W/\opi)$, and the condition that the expression under the
radical in \eqref{a3_03} is positive:
\be
 0 \,<\, \beta W \,<\, \gamma A^2\,.
 \label{a3_04}
 \ee
The maximum value of the growth rate occurs at the midpoint of this interval and is
\be
\lambda_{\rm max}= \gamma A^2.
\label{a3_05}
\ee

Note that the periodicity condition for $p$ and $m$ does {\em not} play here a critical role in
determining the instability increment, in stark contrast to the case of the soliton background
considered in Section 3.2. Namely, as follows from \eqref{e3_08} and \eqref{a3_01}, here this
condition simply requires that the frequency $\omega=\opi-\epsilon W$ fall onto the frequency grid:
 \ $\opi-\epsilon W=2\pi\ell/T$, where $\ell$ is an integer. Thus, if the width of the instability
 band is less than the frequency grid spacing:
 \be
 \gamma A^2/(|\beta|\opi) \,<\, 2\pi/T,
 \label{a3_06}
 \ee
 it is possible that $2\pi\ell/T$ may fall outside the instability band. 
 In this case, the instability will not occur even if $\dz$ exceeds the threshold \eqref{e1_07}.
 This was originally pointed out by Weideman and Herbst \cite{WH} and studied in detail by Yang in \cite{J}.

% -------------------------------------------------------------------------------------

\end{document}